\documentclass[12pt,a4paper]{article}

\usepackage[margin=0.8in]{geometry}
\usepackage{amsmath}
\usepackage{amsmath,bm}
\usepackage{amsfonts}
\usepackage{amssymb}
\usepackage{graphicx}
\usepackage{subfigure}
\usepackage{float}
\usepackage{epstopdf}
\usepackage{cite}
\usepackage{color}
\usepackage{fancyhdr}
\usepackage{authblk}
\usepackage{pdfpages}
\usepackage{pdfpages}
\usepackage{lineno}
\usepackage{csquotes}

\title{Galerkin--Ivanov transformation for nonsmooth modeling of vibro-impacts in continuous structures}
\author[1]{Surya Samukham}
\author[1]{S. N. Khaderi} 
\author[1]{C. P. Vyasarayani \thanks{Address all correspondence to this author. E-mail address: vcprakash@iith.ac.in}}
\affil[1]{Department of Mechanical and Aerospace Engineering, Indian Institute of Technology Hyderabad,
		Kandi, Sangareddy, Telangana, India, 502285}

\begin{document}
	\date{}
	\maketitle
\begin{abstract}
This work deals with the modeling of nonsmooth vibro-impact motion of a continuous structure against a rigid distributed obstacle. Galerkin's approach is used to approximate the solutions of the governing partial differential equations of the structure, which results in a system of ordinary differential equations (ODEs). 
When these ODEs are subjected to unilateral constraints and velocity jump conditions, one must use an event detection algorithm to calculate the time of impact accurately. Event detection in the presence of multiple simultaneous impacts is a computationally demanding task. 
Ivanov (Ivanov, A., $1993$. ``Analytical methods in the theory of vibro-impact systems", Journal of Applied Mathematics and Mechanics, $57(2)$, pp. $221$--$236$.) proposed a nonsmooth transformation for a vibro-impacting multi-degree-of-freedom system subjected to a single unilateral constraint. This transformation eliminates the unilateral constraints from the problem and, therefore, no event detection is required during numerical integration. This nonsmooth transformation leads to sign function nonlinearities in the equations of motion. However, they can be easily accounted for during numerical integration. Ivanov used his transformation to make analytical calculations for the stability and bifurcations of vibro-impacting motions; however, he did not explore its application for simulating distributed collisions in spatially continuous structures. We adopt Ivanov's transformation to deal with multiple unilateral constraints in spatially continuous structures. 
Also, imposing the velocity jump conditions exactly in the modal coordinates is nontrivial and challenging. Therefore, in this work we use a modal-physical transformation to convert the system from modal to physical coordinates on a spatially discretized grid. We then apply Ivanov's transformation on the physical system to simulate the vibro-impact motion of the structure.
The developed method is demonstrated by modeling the distributed collision of a nonlinear string against a rigid distributed surface. For validation, we compare our results with the well-known penalty approach.
\end{abstract}

\maketitle
\section{Introduction}
In many engineering applications, structures are subjected to vibro-impacting motions~\cite{babitsky2013theory,ibrahim2009vibro,ibrahim2014recent}. Clearances in mechanical joints due to wear can lead to vibro-impacting motions in machine components. In some vibration-based energy harvesting applications, vibro-impacting motions are deliberately introduced to increase the bandwidth of the frequency response~\cite{vijayan2015non}. 
In addition, the study of vibro-impact motion of a structure has various applications in the synthesis of sound in stringed musical instruments of Indian origin, like the sitar and the tanpura~\cite{chatziioannou2015energy,issanchou2017modal,issanchou2018nonsmooth}. Applications are also found in the elevator-rope collision problem~\cite{shin2018vibration}, multi-body dynamics~\cite{shafei2018oblique}, and impact-damper systems~\cite{li2006experiments,yang2019investigation}. However, in this work, we are interested in the problem of modeling the vibro-impact motion of a continuous structure constrained by a rigid distributed obstacle. When a continuous structure is spatially discretized using numerical approximation methods, there are three ways the impacting motion can be simulated. 
The first method approximates the rigid obstacle as a foundation of springs with high stiffness and is known as the penalty approach~\cite{gilardi2002literature,bilbao2015numerical,chatziioannou2015energy,bilbao2015numerical2,issanchou2018string,kang2018calculation}. 
Alternatively, one can use Lagrange multipliers to impose contact constraints once the structure comes into contact with the obstacle. Complementarity conditions between Lagrange multipliers and gap functions are used to impose contact loss~\cite{wagg2005periodic}. The Lagrange multiplier approach simulates the perfect sticking motion of vibro-impacting systems~\cite{wagg2005periodic}. Impacts can also be simulated using a coefficient of restitution (CoR or $R$) approach~\cite{gilardi2002literature,li2015modeling,uchida2015making,liakou2016fast,huang2016multi,rebouccas2019unilateral}. Once contact is detected between the structure and the obstacle, the appropriate velocity jump conditions are imposed at the point of contact. In the Lagrange multiplier and CoR approaches, one must solve an event detection problem to isolate the time of impact accurately. Multiple simultaneous impacts are present when a continuous structure impacts a distributed obstacle, which makes the event detection problem computationally expensive. A few other recently developed numerical time-integration techniques to simulate the nonsmooth dynamical systems are given in \cite{acary2008numerical}, \cite{studer2009numerics}, and \cite{acary2016energy}. 

Within the framework of the CoR-based approach, Ivanov~\cite{ivanov1993analytical,ivanov1994impact} proposed a nonsmooth spatial transformation that automatically satisfies both the unilateral constraints and the velocity jump conditions at the point of contact. This approach eliminates the need for event detection. Ivanov proposed this method to study the vibro-impacting motion of a multi-degree-of-freedom (MDOF) system, where the displacement of a single mass is constrained. Ivanov's transformation was successfully used to study the nonlinear dynamics of a single-degree-of-freedom (SDOF) vibro-impacting ship motion~\cite{grace2011inelastic,grace2011inelastic2}. In this work, we adopt Ivanov's transformation to account for the multiple unilateral constraints in a MDOF system. 

In this paper, we have implemented Ivanov's approach to simulate the vibro-impacting motion of a continuous structure against a rigid distributed obstacle. First, we discuss the implementation of Ivanov's transformation to a SDOF system, followed by its application to the MDOF system (continuous structure). To demonstrate the application of Ivanov's approach to a continuous structure subjected to distributed contact, we consider a nonlinear string and present the vibro-impacting motion of the string against the rigid flat and sinusoidal obstacles for different CoR. Galerkin approximation has been used to solve the governing partial differential equations (PDEs) of the string. The set of ordinary differential equations (ODEs) obtained from the Galerkin approach are in modal coordinates and the imposition of impact constraints upon the modal system is a challenging task~\cite{vyasarayani2010modeling,wagg2002application,wei2019effect}. To avoid such difficulties, we use a transformation to convert the modal system into its physical coordinates by discretizing the string in space \cite{van2016simulation,issanchou2017modal}. 
We then apply Ivanov's transformation to the spatially discretized system to incorporate the impact constraints.
After Ivanov's transformation, the resulting system of differential equations are in Ivanov's coordinates and non-stiff in nature. 
Therefore, an explicit integrator with a large time-step size can be used to solve the system of differential equations, which is not the case with the penalty method because of the presence of stiff differential equations.
This advantage of implementing an explicit integrator with a larger time-step size in a MDOF system significantly reduces the computational cost and time compared to an implicit integrator.
The results obtained from Ivanov's method in this work have been verified with the penalty method for the case of $R=1$.
The main advantages of the proposed method are as follows:
\begin{enumerate}
\item The nonsmooth nature of the contact problem is preserved, which is not the case with the penalty approach.
\item The equations in Ivanov's coordinates are non-stiff, in contrast to those obtained using the penalty approach.
\item Explicit integration schemes can be used with larger time-step sizes for integration purposes instead of a computationally expensive implicit integrator.
\item No event detection is necessary for simulating impacts, which is the most important benefit of Ivanov's approach, since the problem of event detection is challenging in the case of large MDOF systems (like Galerkin approximations and finite element analysis of continuous structures) subjected to distributed collisions.
\end{enumerate}

This paper is organized as follows: In Sec.~{\ref{Model}}, we describe in detail how Ivanov's transformation works by modeling a point mass falling on a rigid obstacle. Using the example of a vibro-impacting SDOF system, we compare the advantages of Ivanov's method over the penalty method in Sec.~\ref{vibroimp}. The equations of motion for a nonlinear MDOF system in Ivanov's coordinates are derived in Sec.~\ref{mainder}. In Sec.~\ref{galerkin}, we discuss the Galerkin approximation of the governing differential equation of a nonlinear string and the modal-to-physical coordinate transformation. In Sec.~\ref{results}, we apply Ivanov's transformation to simulate the vibro-impact motion of the string against a distributed obstacle. The results are validated and compared with the penalty approach. The contributions of the paper have been summarized in Sec.~\ref{concl}.

\section{Mathematical modeling}
In this section, using the example of a point mass bouncing on a rigid surface, we illustrate the idea behind Ivanov's transformation. Later, by using the example of a vibro-impacting oscillator, we demonstrate the advantages of modeling impact using Ivanov's method over the penalty approach. When simulating rigid collisions using a penalty approach, the penalty stiffness term is usually selected to be large, which results in stiff differential equations. Due to the stiff nature of the differential equations, one must use a very small time-step size during numerical integration. However, in Ivanov's method, the obtained differential equations in the transformed coordinates are non-stiff. Unlike the penalty approach, Ivanov's method accurately captures the nonsmooth behavior in velocity during impact. 

\subsection{Point mass bouncing on a rigid surface}
\label{Model}
\begin{figure}[htpb!]
\begin{center}
\includegraphics[width=0.40\textwidth]{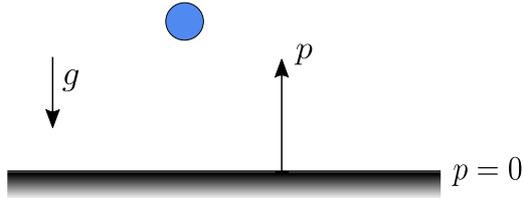}
\caption{Schematic of a point mass falling on a rigid surface.}
\label{Figure1}
\end{center}
\end{figure}

\begin{figure}[htpb!]
\begin{center}
\subfigure[]{\includegraphics[width=0.30\textwidth]{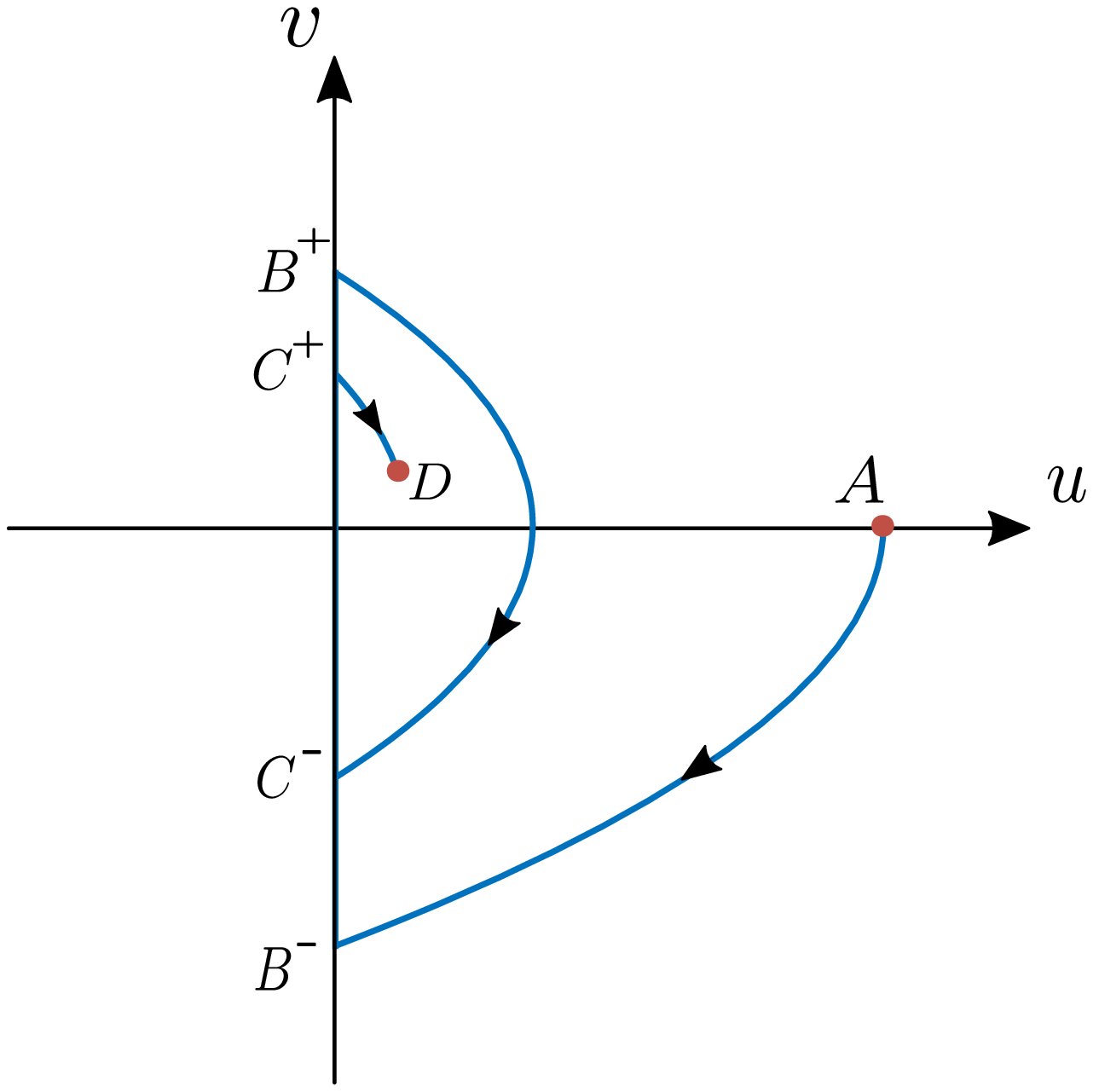}}
\subfigure[]{\includegraphics[width=0.40\textwidth]{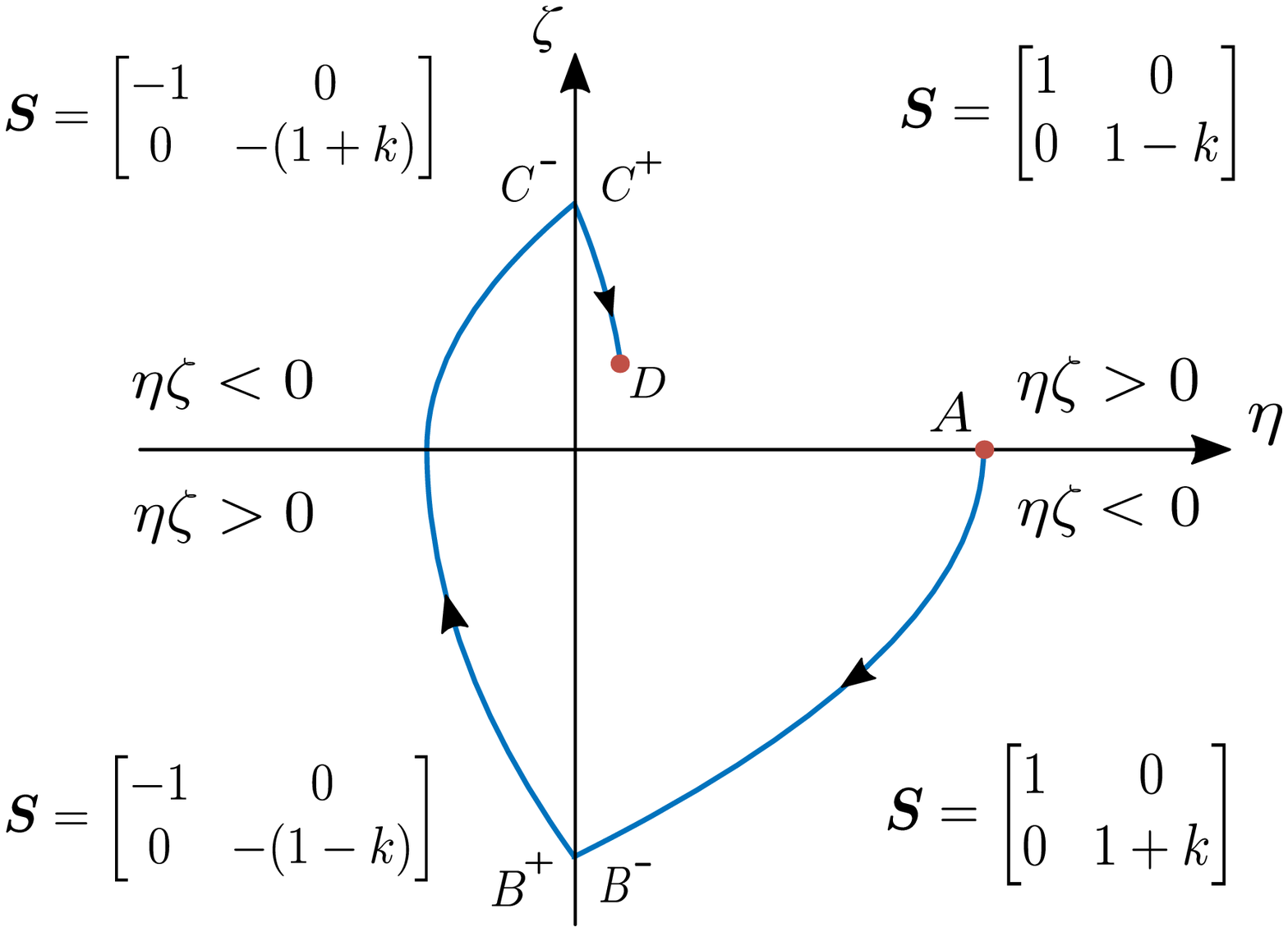}}
\caption{(a) Phase space in original $u-v$ coordinates. (b) Phase space in transformed $\eta-\zeta$ coordinates.}
\label{Figure2}
\end{center}
\end{figure}

\begin{figure}[htpb!]
\centering
\subfigure[]{\includegraphics[width=0.45\textwidth]{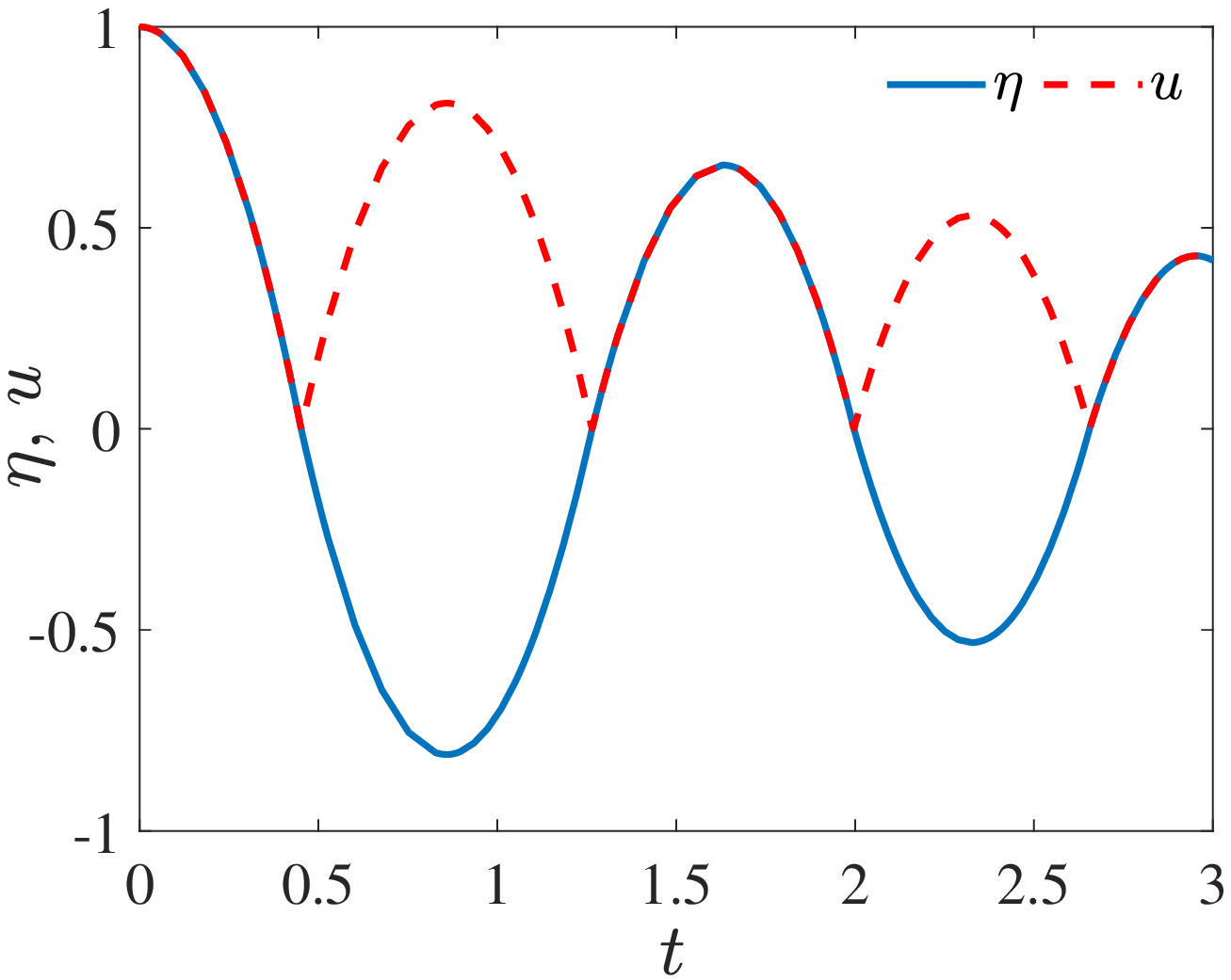}\label{Figure3a}}
\subfigure[]{\includegraphics[width=0.45\textwidth]{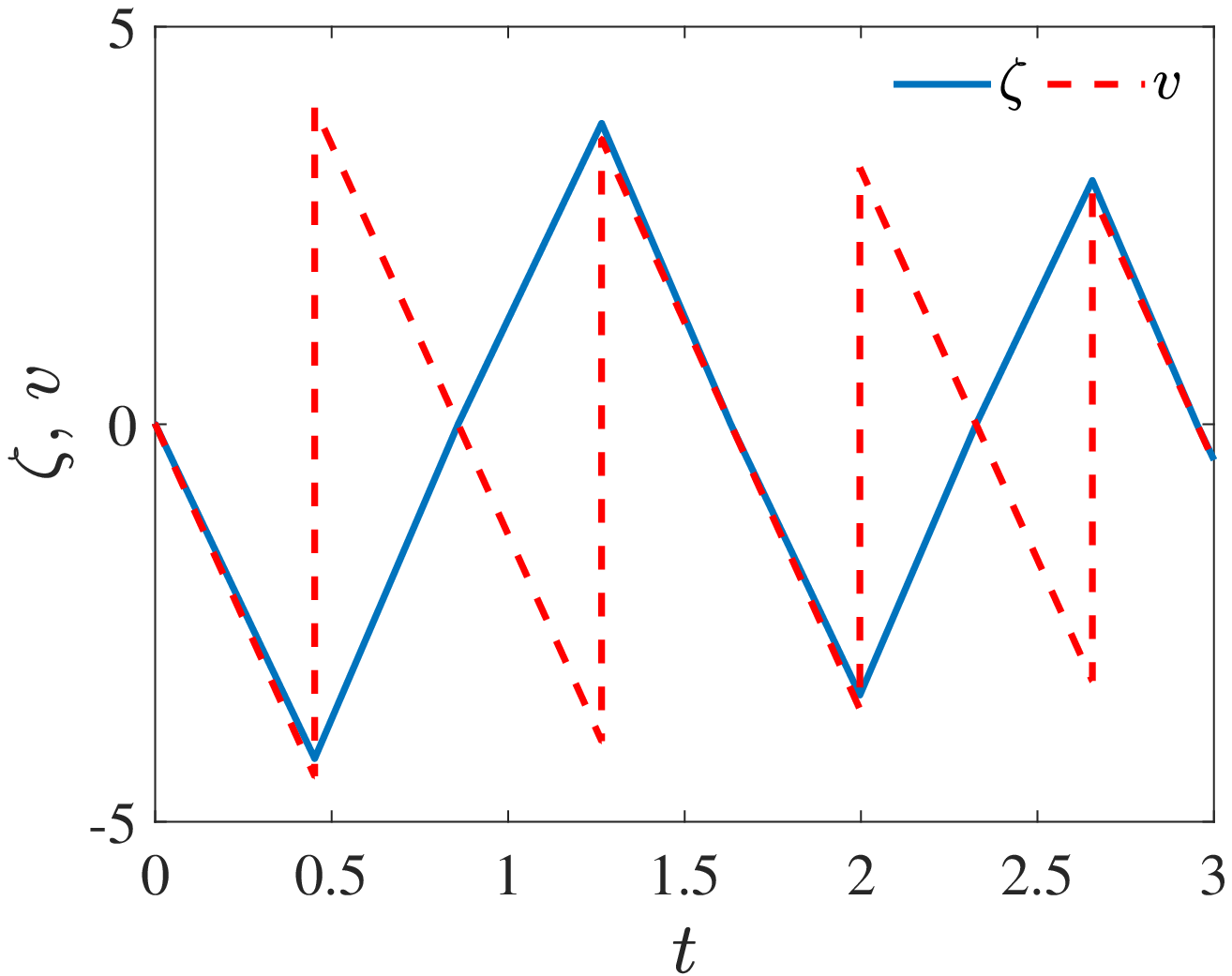}\label{Figure3b}}
\subfigure[]{\includegraphics[width=0.45\textwidth]{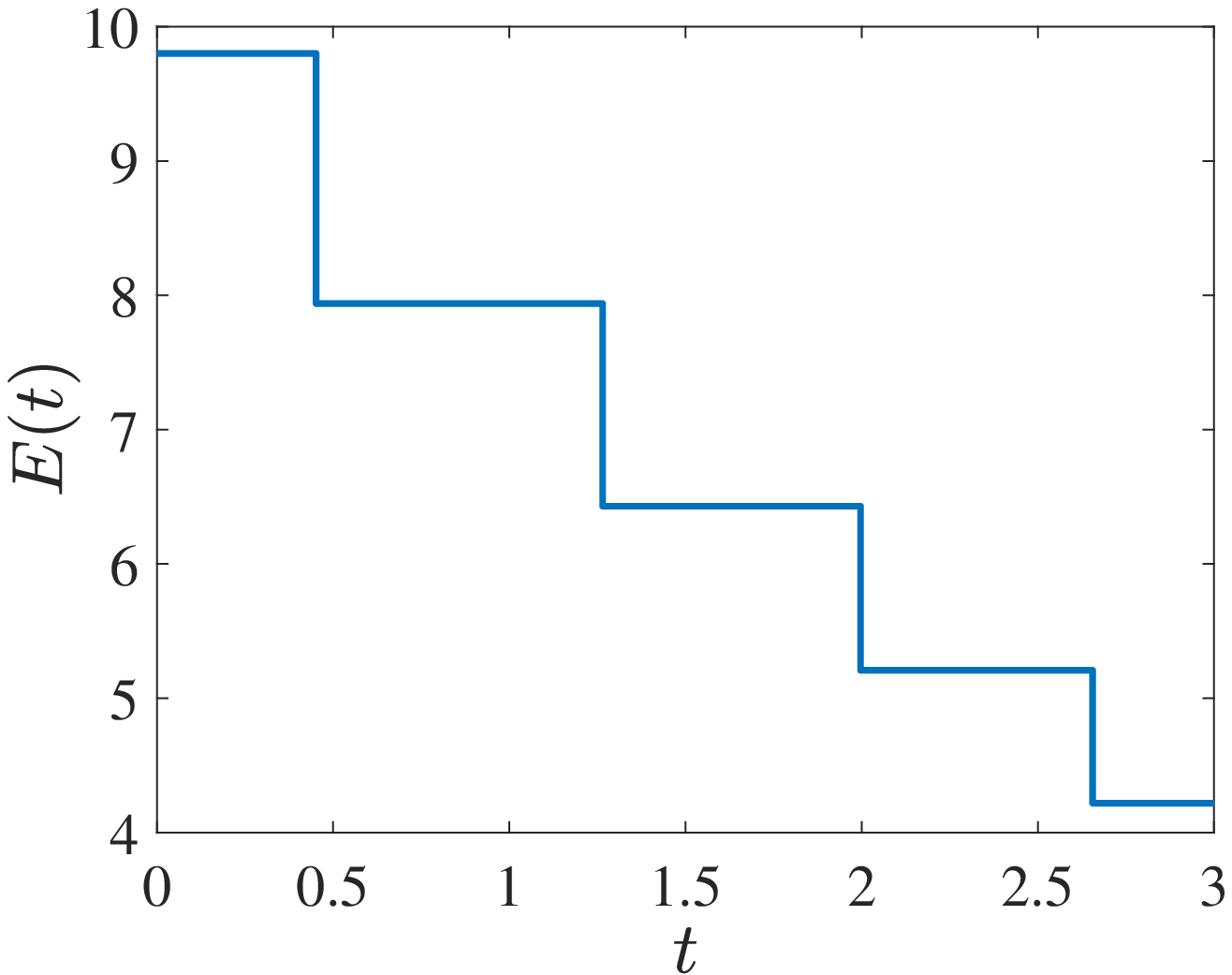}\label{Figure3c}}
\caption{(a) Position of the point mass in physical coordinates ($u$) and transformed coordinates ($\eta$). (b) The velocity of the point mass in physical ($v$) and transformed  ($\zeta$) coordinates. (c) Energy of the point mass over time. The results were obtained for $u(0)=1$, $v(0)=0$, $g=9.8$, and $R=0.9$.}
\label{Pointmassfig}
\end{figure}

To illustrate Ivanov's transformation, we consider the problem of a point mass, under the influence of gravity, impacting a rigid surface (see Fig.~\ref{Figure1}). The equation of motion for this problem can be written as follows:
\begin{equation}
\ddot{p}=-g, \thinspace\thinspace\thinspace p\ge0,
\label{Sec.2.1.eq1}
\end{equation}
with initial conditions $p(0)=\alpha_0$ and $\dot{p}(0)=\beta_0$. When the mass makes contact with the obstacle ($p(t_c)=0$) at time $t_c$, the following jump condition in velocity must be respected:
\begin{equation}
\dot{p}(t_{c}^{+})=-R\dot{p}(t_{c}^{-}),
\label{Sec.2.1.eq2}
\end{equation}
where $R$ is the coefficient of restitution. Equations~(\ref{Sec.2.1.eq1}) and (\ref{Sec.2.1.eq2}) completely describe the motion of the mass. To apply Ivanov's transformation, we recast Eqs.~(\ref{Sec.2.1.eq1}) and (\ref{Sec.2.1.eq2}) in state-space form. Introducing $u=p$ and $v=\dot{p}$, we get
\begin{equation}
\left\{ \begin{array}{c}
\dot{u}\\ 
\dot{v}
\end{array} \right\} =\left[\begin{array}{cc}
0 & 1\\
0 & 0
\end{array}\right]\left\{ \begin{array}{c}
u\\
v
\end{array}\right\} +\left\{ \begin{array}{c}
0\\
-g
\end{array}\right\}, \thinspace\thinspace\thinspace u\ge0
\label{Sec.2.1.eq3}
\end{equation}
and when $u(t_c)=0$, we have
\begin{equation}
v(t_{c}^{+})=-Rv(t_{c}^{-})
\label{Sec.2.1.eq4}
\end{equation}
Ivanov introduced the following nonsmooth transformation:
\begin{eqnarray}
\label{Sec.2.1.eq5}
&u =\eta ~\text{sgn}(\eta)\\
\label{Sec.2.1.eq6}
&v =\zeta\left(1-k~ \text{sgn}(\eta\zeta)\right)\text{sgn}(\eta)
\end{eqnarray}
In Eq.~(\ref{Sec.2.1.eq6}), $k=\frac{1-R}{1+R}$. Equation~(\ref{Sec.2.1.eq5}) ensures that $u\ge0$ for all time. To understand the transformation given by Eq.~(\ref{Sec.2.1.eq6}), we need to look at the phase-space of $\eta-\zeta$, as shown in  Fig.~\ref{Figure2}(b); the actual phase-space (in $u-v$) is shown in Fig.~\ref{Figure2}(a) for comparison. In the first and third quadrants of the $\eta-\zeta$ plane, we have $\text{sgn}(\eta\zeta)>0$. In the second and fourth quadrants of the $\eta-\zeta$ plane, we have $\text{sgn}(\eta\zeta)<0$. When a trajectory $AB^{-}$ starting in phase-space $\eta-\zeta$ from point $A$  is about to cross $\eta=0$ at time $t_{B}^{-}$ from the fourth to the third quadrant, from Eq.~(\ref{Sec.2.1.eq6}), we have
\begin{equation}
\label{Sec.2.1.eq7}
v(t_{B}^{-})=\left(1+k\right)\zeta(t_{B}^{-})
\end{equation}
Immediately after the contact, at time $t_B^{+}$, we have
\begin{equation}
\label{Sec.2.1.eq8}
v(t_{B}^{+})=-\left(1-k\right)\zeta(t_{B}^{+})
\end{equation}
Similarly, when the trajectory $B^{+}C^{-}$ is about to cross $\eta=0$ at time $t_{C}^{-}$ from the second to the first quadrant, we have
\begin{equation}
\label{Sec.2.1.eq9}
v(t_{C}^{-})=-\left(1+k\right)\zeta(t_{C}^{-})\\
\end{equation}
After crossing $\eta=0$, we have
\begin{equation}
\label{Sec.2.1.eq10}
v(t_{C}^{+})=\left(1-k\right)\zeta(t_{C}^{+})
\end{equation}
It should be noted that both $\eta$ and $\zeta$ are continuous, so we have $\zeta(t_{B}^{-})=\zeta(t_{B}^{+})$ and $\zeta(t_{C}^{-})=\zeta(t_{C}^{+})$, which results in the following expressions:
\begin{eqnarray}
\label{Sec.2.1.eq111}
v(t_{B}^{+})=-\left(\frac{1-k}{1+k}\right)v(t_{B}^{-})=-Rv(t_{B}^{-})\\
\label{Sec.2.1.eq112}
v(t_{C}^{+})=-\left(\frac{1-k}{1+k}\right)v(t_{C}^{-})=-Rv(t_{C}^{-})
\end{eqnarray}
Therefore, it is clear that the nonsmooth transformation given by Eq.~(\ref{Sec.2.1.eq6}) automatically satisfies Eq.~(\ref{Sec.2.1.eq4}) at the event of an impact. Due to the jump conditions (Eq.~(\ref{Sec.2.1.eq111}) and Eq.~(\ref{Sec.2.1.eq112})) imposed on the velocity $v$ at the time of impact, the trajectories in $\eta-\zeta$ space that enter the third quadrant from the fourth quadrant can only go to the first quadrant through the second quadrant. 
Equations~(\ref{Sec.2.1.eq5}) and~(\ref{Sec.2.1.eq6}) can be written as follows:
\begin{equation}
\left\{ \begin{array}{c}
u\\
v
\end{array}\right\} =\boldsymbol{S}\left\{ \begin{array}{c}
\eta\\
\zeta
\end{array}\right\}
\label{Sec.2.1.eq12}
\end{equation}
The transformation matrix, $\boldsymbol{S}$ in Eq.~(\ref{Sec.2.1.eq12}), is given by
\begin{equation}
\label{Sec.2.1.eq13}
\boldsymbol{S}=\left[\begin{array}{cc}
\text{sgn}(\eta) & 0\\
0 & \left(1-k~\text{sgn}(\eta \zeta)\right)\text{sgn}(\eta)
\end{array}\right]
\end{equation}
The matrix $\boldsymbol{S}$ is constant in each quadrant of phase-space (see Fig.~\ref{Figure2}). Substituting Eq.~(\ref{Sec.2.1.eq12}) into Eq.~(\ref{Sec.2.1.eq3}), we get
\begin{equation}
\left\{ \begin{array}{c}
\dot{\eta}\\
\dot{\zeta}
\end{array}\right\} =\boldsymbol{S}^{-1}\left[\begin{array}{cc}
0 & 1\\
0 & 0
\end{array}\right]\boldsymbol{S}\left\{ \begin{array}{c}
\eta\\
\zeta
\end{array}\right\} +\boldsymbol{S}^{-1}\left\{ \begin{array}{c}
0\\
-g
\end{array}\right\} 
\label{Sec.2.1.eq14}
\end{equation}
The initial conditions for $\eta$ and $\zeta$ can be calculated as follows:
\begin{equation}
\left\{ \begin{array}{c}
\eta(0)\\
\zeta(0)
\end{array}\right\} =\boldsymbol{S}^{-1}\left\{ \begin{array}{c}
u(0)\\
v(0)
\end{array}\right\} 
\label{Sec.2.1.eq15}
\end{equation}
Equation~(\ref{Sec.2.1.eq14}) can be integrated numerically to obtain $\eta$ and $\zeta$. Later, Eq.~(\ref{Sec.2.1.eq12}) can be used to calculate $u$ and $v$. It should be noted that there are no constraints on Eq.~(\ref{Sec.2.1.eq14}) and its solutions are continuous functions of time. Figure~\ref{Figure3a} shows the displacement of the point mass in transformed coordinates ($\eta$, solid blue line) and physical coordinates ($u$, dashed red line). Similarly, Fig.~\ref{Figure3b} shows the velocity in transformed coordinates ($\zeta$, solid blue line) and physical ($v$, dashed red line) coordinates. Figure~\ref{Figure3c} shows the total energy of the point mass, and it can be seen that energy is lost at every collision. If we integrate Eq.~(\ref{Sec.2.1.eq14}) for a sufficiently long time, the energy will approach zero. Also, the frequency of the transformed solution ($\eta$) continuously increases due to the chattering nature of the physical solution ($u$). If we use an adaptive time-step integrator, the integrator automatically reduces the time-step sizes to account for increasing frequency in the solution. For practical reasons, if the energy falls below a small threshold, we can assume that the point mass has reached equilibrium.

\subsection{Single-degree-of-freedom vibro-impacting system}
\label{vibroimp}
\begin{figure}[htpb!]
\begin{center}
\subfigure[]{\includegraphics[width=0.39\textwidth]{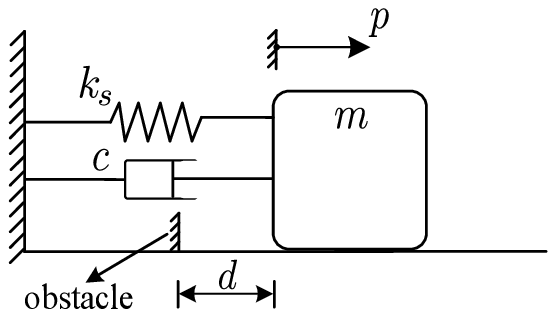}\label{Figure4a}} \hspace{2cm}
\subfigure[]{\includegraphics[width=0.39\textwidth]{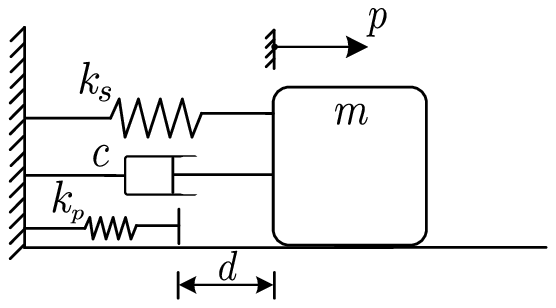}\label{Figure4b}}
\caption{(a) Schematic of a spring-mass system with an obstacle. (b) Equivalent representation of the spring-mass system with the rigid obstacle replaced by a linear spring.}
\label{Figure4}
\end{center}
\end{figure}

To demonstrate the advantage of modeling impact using Ivanov's method over a penalty approach, we consider the problem of a SDOF vibro-impacting system, as shown in Fig.~\ref{Figure4}(a). In the penalty approach, the rigid obstacle as shown in  Fig.~\ref{Figure4}(a) is replaced with a high-stiffness ($k_p$) linear spring (see Fig.~\ref{Figure4}(b)), usually referred to as the penalty stiffness. For a large value of $k_p$, the spring behaves approximately like a rigid obstacle. In a penalty approach, we solve the following equation:
\begin{equation}
\label{Sec.2.2.eq1}
m\ddot{p}+c\dot{p}+(k_{s}+\mu k_{p})p=0
\end{equation}
In Eq.~(\ref{Sec.2.2.eq1}), $\mu$ is defined as follows:
\begin{equation}
\mu=\begin{cases}
0 & p\ge 0\\
1 & p<0
\end{cases}
\label{Sec.2.2.eq2}
\end{equation}
From Eq.~(\ref{Sec.2.2.eq1}), we can see that the system has a time-scale of $\tau_1=2\pi \sqrt{\frac{m}{k_s}}$ when the mass is not in contact with the obstacle ($\mu=0$); when the mass is in contact with the obstacle ($\mu=1$), it has a time-scale of $\tau_2=2\pi \sqrt{\frac{m}{k_s+k_{p}}}$. Having two separate time-scales of very different magnitudes makes Eq.~(\ref{Sec.2.2.eq1}) a stiff differential equation. Therefore, when an explicit integrator is used to solve Eq.~(\ref{Sec.2.2.eq1}), one must use a time-step size much smaller than $\tau_2$ (the smallest time-scale present in the solution) to obtain accurate results. 
The equation of motion of the system, considering a perfectly rigid obstacle ($R=1$) (see Fig.~\ref{Figure4a}), can be written as follows:
\begin{equation}
\label{Sec.2.2.eq3}
m\ddot{p}+c\dot{p}+k_{s}p=0, ~~p\ge d
\end{equation}
To apply Ivanov's method for solving the vibro-impact problem (Eq.~(\ref{Sec.2.2.eq3})), we substitute $u=p-d$ and $v=\dot{p}$ in Eq.~(\ref{Sec.2.2.eq3}) to obtain
\begin{eqnarray}
\label{Sec.2.2.eq4}
\left\{ \begin{array}{c}
\dot{u}\\
\dot{v}
\end{array}\right\} =\left[\begin{array}{cc}
0 & 1\\
-k_s/m & -c/m
\end{array}\right]\left\{ \begin{array}{c}
u\\
v
\end{array}\right\} \nonumber \\ + \left\{ \begin{array}{c}
0\\
-k_s d/m 
\end{array}\right\},\,~u\ge0 
\end{eqnarray}
By following the same transformation as discussed in the point mass bouncing problem (see Eqs.~(\ref{Sec.2.1.eq5}) and (\ref{Sec.2.1.eq6})), we can write the equations of motion for the vibro-impacting system with perfectly rigid impacts in Ivanov's coordinates as follows:
\begin{eqnarray}
\left\{ \begin{array}{c}
\dot{\eta}\\
\dot{\zeta}
\end{array}\right\} =\boldsymbol{S}^{-1}\left[\begin{array}{cc}
0 & 1\\
-k_s/m & -c/m
\end{array}\right]\boldsymbol{S}\left\{ \begin{array}{c}
\eta\\
\zeta
\end{array}\right\} \nonumber \\ + \boldsymbol{S}^{-1} \left\{ \begin{array}{c}
0\\
-k_s d/m
\end{array}\right\}
\label{Sec.2.2.eq5}
\end{eqnarray}
(Refer to Eq.~(\ref{Sec.2.1.eq13}) for the definition of the matrix $\boldsymbol{S}$.) As pointed out earlier (see Fig.~\ref{Figure2}), the matrix $\boldsymbol{S}$ is diagonal and constant in each quadrant of the phase-space. Further, the magnitude of the elements of $\boldsymbol{S}$ are of the same order in all quadrants (see Fig.~\ref{Figure2}). Therefore, when Eq.~(\ref{Sec.2.2.eq5}) is solved using an explicit integrator, it does not appear as a stiff equation to the integrator. 

\begin{figure}[htpb!]
\centering
\includegraphics[width=0.5\textwidth]{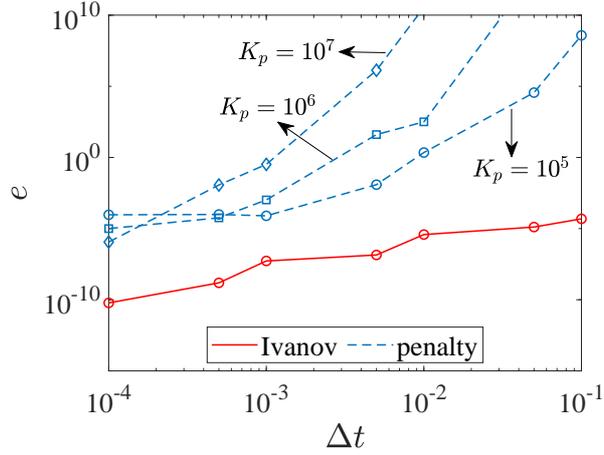}
\caption{Variation of mean squared error ($e=\frac{1}{n}\sum_{i=1}^{n}(p_r(t_i)-p(t_i))^2$) with respect to the time-step size of integration ($\Delta t$) in Ivanov's method and the penalty method.}
\label{Figure5}
\end{figure}

First, we determine the solution of the system shown in Fig.~\ref{Figure4a} subject to the constraint $p \ge 0.5$ using the event detection technique. The solution of Eq.~(\ref{Sec.2.2.eq3}) is obtained using the ``ode45" (adaptive time-step) integrator in MATLAB with relative and absolute tolerances of $1\times10^{-12}$. The parameters considered for the analysis are $m=1$, $k=1$, $c=0$, and the integration has been carried out over a time of $t=0$ to $t=10$. The initial conditions are $p(0)=1$ and $\dot{p}(0)=0$.
This reference solution $p_r(t)$ is used to compare the results from Ivanov's method and the penalty method with an explicit fixed-time-step-size integrator.
Now, we solve the same problem (Fig.~\ref{Figure4}) using the penalty approach (Eq.~(\ref{Sec.2.2.eq1})) and Ivanov's method (Eq.~(\ref{Sec.2.2.eq5})) using a fixed-time-step-size, fourth-order Runge--Kutta integrator. Different time-step sizes ($\Delta t$) in the range of $10^{-4}$ to $10^{-1}$ are considered for the analysis. Further, the mean squared error (MSE) given by $e=\frac{1}{n}\sum_{i=1}^{n}(p_r(t_i)-p(t_i))^2$ is evaluated between the reference solution $p_r(t)$ and the solutions obtained from Ivanov's method and the penalty method ($p(t)$) for different time-step sizes of integration. Here, $n$ is the number of sample points considered for the error estimation and $t_i$ is the time instant. The variation of MSE for the two methods with respect to the change in integration time-step size is shown in Fig.~\ref{Figure5}. In Fig.~\ref{Figure5}, $e$ is calculated by evaluating the solutions ($p(t)$ and $p_r(t)$) at $n=101$ equally spaced time intervals between $t=0$ and $t=10$.
The MSE for the penalty approach is calculated for $3$ different penalty stiffness values ($k_p= 10^5, 10^6,$ and $10^7$) and the corresponding variation in $e$ is shown in Fig.~\ref{Figure5}.
From the results shown in Fig.~\ref{Figure5}, it is clearly observed that the MSE in the penalty method increases with an increase in $\Delta t$.
This is because of the penalty stiffness $k_p$ in Eq.~(\ref{Sec.2.2.eq1}), which decreases the fastest time-scale $\tau_2$ and makes the equation stiff.
Therefore, as discussed earlier, when an explicit integrator is used for the integration, $\Delta t$ must be chosen much smaller than $\tau_2$.
As a result, in the penalty method, the suitable time-step size for integration decreases with an increase in the penalty stiffness $k_p$, which can be clearly observed from Fig.~\ref{Figure5}.
Also, with the increase in $k_p$, the MSE in the solution decreases for the smaller integration time-step sizes (Fig.~\ref{Figure5}). Therefore, to obtain an accurate solution in the penalty method, a very large value of $k_p$ must be chosen. 
However, this decreases $\tau_2$ and a correspondingly smaller time-step size must be used for integration, which increases the computation cost significantly.
In contrast, the final equations obtained in Ivanov's method are not stiff (see Eq.~(\ref{Sec.2.2.eq5})) and hence a larger time-step size can be used for integration without reducing the accuracy of the solution. We can clearly see this in Fig.~\ref{Figure5}: the error in the solution obtained from Ivanov's method is of an order smaller than $10^{-4}$ even when $\Delta t=10^{-1}$, which is not the case with the penalty approach.

In summary, with the penalty approach, to obtain an accurate solution of a vibro-impact system, the time-step size must chosen to be very small depending on the penalty stiffness $k_p$. However, by using Ivanov's method, an accurate solution can be obtained with a much larger time-step size. Moreover, for simulating a rigid collision, one must use a large value for the penalty stiffness, which puts a hard restriction on the numerical time-step size of the integrator. However, in Ivanov's approach, the nonsmooth transformation simulates the case of infinite penalty stiffness without making the equations stiff. 

\subsection{Ivanov's transformation for distributed impacts}
\label{mainder}
Having demonstrated the applicability and advantage of Ivanov's transformation for an SDOF system, we now discuss its extension to MDOF systems where only some degrees of freedom are subjected to impact constraints. 
To develop Ivanov's transformation for the most general MDOF system, we consider the following coupled nonlinear model where some degrees of freedom ($\boldsymbol{p}$) are subjected to impact while the others ($\boldsymbol{q}$) are not subjected to any impact. The equations of motion for such a mechanical system can be written as follows:
\begin{eqnarray}
\label{Sec.2.3.eq1}
\ddot{\boldsymbol{p}}= & \boldsymbol{f}(\boldsymbol{p},\dot{\boldsymbol{p}},\boldsymbol{q},\dot{\boldsymbol{q}},t)\\
\label{Sec.2.3.eq2}
\ddot{\boldsymbol{q}}= & \boldsymbol{g}(\boldsymbol{p},\dot{\boldsymbol{p}},\boldsymbol{q},\dot{\boldsymbol{q}},t),
\end{eqnarray}
with the initial conditions $\boldsymbol{p}(\boldsymbol{0})=\boldsymbol{\alpha}_0$, $\dot{\boldsymbol{p}}(\boldsymbol{0})=\boldsymbol{\beta}_0$, $\boldsymbol{q}(\boldsymbol{0})=\boldsymbol{\gamma}_0$, and $\dot{\boldsymbol{q}}(\boldsymbol{0})=\boldsymbol{\nu}_0$. In Eq.~(\ref{Sec.2.3.eq1}), the displacement coordinates $\boldsymbol{p}=[p_{1}(t),\thinspace p_{2}(t),...,p_{m}(t)]^\text{T}$ are subjected to impact constraints of the form
\begin{equation}
p_{i}\ge d_{i},\thinspace\thinspace i=1,2,...,m
\label{Sec.2.3.eq3}
\end{equation}
Once any of the above constraints becomes an equality constraint ($p_{i}(t_c)= d_{i}$) at time $t_c$, a velocity jump condition is imposed as follows:
\begin{equation}
\dot{p}_{i}(t_c^{+})=-R\dot{p}_{i}(t_c^{-})
\label{Sec.2.3.eq4}
\end{equation}
By defining $\boldsymbol{d}=[d_1,d_2,...,d_m]^\text{T}$, Eq.~(\ref{Sec.2.3.eq3}) can be written compactly as
\begin{equation}
\boldsymbol{p}-\boldsymbol{d} \ge \boldsymbol{0}
\label{Sec.2.3.eq5}
\end{equation}
It should be noted that the displacement coordinates $\boldsymbol{q}=[q_{1}(t),\thinspace q_{2}(t),...,q_{n}(t)]^\text{T}$ are not subjected to any impact constraints. By introducing the state variables  $\boldsymbol{u}=\boldsymbol{p}-\boldsymbol{d}$, $\boldsymbol{v}=\dot{\boldsymbol{p}}$, $\boldsymbol{r}=\boldsymbol{q}$, and $\boldsymbol{s}=\dot{\boldsymbol{q}}$, Eqs.~(\ref{Sec.2.3.eq1}) and~(\ref{Sec.2.3.eq2}) can be written as follows:
\begin{eqnarray}
\label{Sec.2.3.eq6}
\dot{\boldsymbol{u}}&=&\boldsymbol{v}\\ \label{eq2_vcp}
\dot{\boldsymbol{v}}&=&\boldsymbol{f}(\boldsymbol{u}+\boldsymbol{d},\boldsymbol{v},\boldsymbol{r},\boldsymbol{s},t)\\ \label{Sec.2.3.eq7}
\dot{\boldsymbol{r}}&=&\boldsymbol{s}\\ \label{eq4_vcp}
\dot{\boldsymbol{s}}&=&\boldsymbol{g}(\boldsymbol{u}+\boldsymbol{d},\boldsymbol{v},\boldsymbol{r},\boldsymbol{s},t)
\end{eqnarray}
The impact constraint becomes $\boldsymbol{u} \ge 0$, where $\boldsymbol{u}=[u_{1}(t),\thinspace u_{2}(t),...,u_{m}(t)]^\text{T}$. Once any of the above constraints becomes an equality constraint ($u_{i}(t_c)=0$) at time $t_c$, a velocity jump condition is imposed as follows:
\begin{equation}
\label{Sec.2.3.eq8}
v_{i}(t_c^{+})=-Rv_{i}(t_c^{-})
\end{equation}
To apply Ivanov's method, the following transformation is introduced:
\begin{eqnarray}
\label{Sec.2.3.eq9}
\boldsymbol{u}&=&\boldsymbol{T}\boldsymbol{\eta}\\\label{Sec.2.3.eq10}
\boldsymbol{v}&=&\boldsymbol{W}\boldsymbol{\zeta}
\end{eqnarray}
The matrices $\boldsymbol{T}$ and $\boldsymbol{W}$ are defined as follows:
\begin{eqnarray}
\label{Sec.2.3.eq11}
\boldsymbol{T}&=&\text{diag}\left(\text{sgn}(\boldsymbol{\eta})\right)\\
\label{Sec.2.3.eq12}
\boldsymbol{W}&=&\text{diag}\left(\left(1-k~\text{sgn}(\boldsymbol{\eta}\circ\boldsymbol{\zeta})\right)\circ\text{sgn}(\boldsymbol{\eta})\right)
\end{eqnarray}
In Eq.~(\ref{Sec.2.3.eq12}), the symbol \enquote{$\circ$} represents the element-by-element multiplication (Hadamard product) of two vectors. Substituting Eqs.~(\ref{Sec.2.3.eq9}) and (\ref{Sec.2.3.eq10}) into Eqs.~(\ref{Sec.2.3.eq6})-(\ref{eq4_vcp}), we get
\begin{eqnarray}
\label{Sec.2.3.eq13}
\dot{\boldsymbol{\eta}}&=&\boldsymbol{T}^{-1}\boldsymbol{W}\boldsymbol{\zeta}\\ \label{Sec.2.3.eq14}
\dot{\boldsymbol{\zeta}}&=&\boldsymbol{W}^{-1}\boldsymbol{f}(\boldsymbol{T}\boldsymbol{\eta}+\boldsymbol{d},\boldsymbol{W}\boldsymbol{\zeta},\boldsymbol{r},\boldsymbol{s},t)\\ \label{Sec.2.3.eq15}
\dot{\boldsymbol{r}}&=&\boldsymbol{s}\\ \label{Sec.2.3.eq16}
\dot{\boldsymbol{s}}&=&\boldsymbol{g}(\boldsymbol{T}\boldsymbol{\eta}+\boldsymbol{d},\boldsymbol{W}\boldsymbol{\zeta},\boldsymbol{r},\boldsymbol{s},t) 
\end{eqnarray}
The initial conditions for Eqs.~(\ref{Sec.2.3.eq13})--(\ref{Sec.2.3.eq16}) are
\begin{eqnarray}
\label{Sec.2.3.eq17}
\boldsymbol{\eta}(0)&=&\boldsymbol{T}^{-1}\boldsymbol{u}(0)=\boldsymbol{T}^{-1}\left(\boldsymbol{p}(0)-\boldsymbol{d}\right) \\ \nonumber
&=&\boldsymbol{T}^{-1}\left(\boldsymbol{\alpha}_{0}-\boldsymbol{d}\right)\\ \label{Sec.2.3.eq18}
\boldsymbol{\zeta}(0)&=&\boldsymbol{W}^{-1}\boldsymbol{v}(0)=\boldsymbol{W}^{-1}\dot{\boldsymbol{p}}(0)=\boldsymbol{W}^{-1}\boldsymbol{\beta}_{0}\\ \label{Sec.2.3.eq19}
\boldsymbol{r}(0)&=&\boldsymbol{q}(0)=\boldsymbol{\gamma}_{0}\\ \label{Sec.2.3.eq20}
\boldsymbol{s}(0)&=&\boldsymbol{\dot{q}}(0)=\boldsymbol{\nu}_{0} 
\end{eqnarray}
The solution of Eqs.~(\ref{Sec.2.3.eq13})--(\ref{Sec.2.3.eq16}) automatically satisfies all the constraints imposed on $\boldsymbol{p}$ (Eq.~(\ref{Sec.2.3.eq3})) and $\boldsymbol{\dot{p}}$ (Eq.~(\ref{Sec.2.3.eq4})). Equations~(\ref{Sec.2.3.eq13})--(\ref{Sec.2.3.eq16}) are the closed-form equations for a nonlinear MDOF system (where only some degrees of freedom are subjected to constraints) in Ivanov's coordinates. Upon solving the system in Ivanov's coordinates, the physical solution (displacements and velocities) can be reconstructed by using the inverse of the transformations given by Eqs.~(\ref{Sec.2.3.eq9}) and (\ref{Sec.2.3.eq10}).

\section{Vibro-impact motion of a string}
\label{galerkin}
\begin{figure}[htpb!]
\begin{center}
\includegraphics[width=0.5\textwidth]{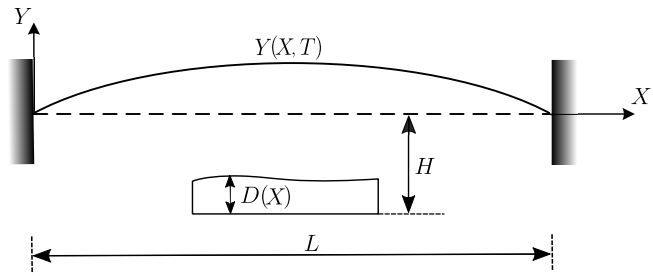}
\caption{Schematic of a string with a distributed impact constraint.}
\label{Figure7}
\end{center}
\end{figure}

To demonstrate the application of the proposed methodology, we simulate the vibro-impact motion of a nonlinear string against a rigid distributed obstacle (see Fig.~\ref{Figure7}). The governing equation of motion for the string is given by~\cite{oplinger1960frequency}
\begin{equation}
\rho \frac{\partial^2 Y}{\partial T^2} - T_0 \left(1+ \frac{EA}{2T_0 L} \int_{0}^{L} \left(\frac{\partial Y}{\partial X}\right)^2 dX\right) \frac{\partial^2 Y}{\partial X^2}+C\frac{\partial Y}{\partial T}=0,
\label{galerkin_eq1}
\end{equation}
with the following boundary conditions:
\begin{equation}
Y(0,T)=Y(L,T)=0
\label{galerkin_eq2}
\end{equation}
Here, $Y$ is the transverse displacement of the string, $L$ is the length, $A$ is its cross-sectional area, $T_0$ is the string tension at the equilibrium position, $\rho$ is the mass density, $C$ is the damping coefficient,  $E$ is the Young's modulus of the material of the string, $X$ is the coordinate along the length of the string, and $T$ is time.

\begin{figure}[htpb!]
\begin{center}
\includegraphics[width=0.45\textwidth]{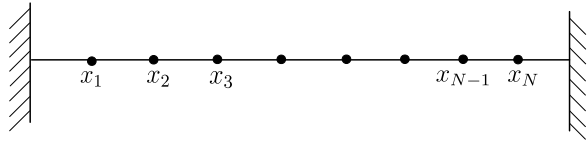}
\caption{Schematic of a spatially discretized string.}
\label{1d_grid}
\end{center}
\end{figure}

We now substitute the following non-dimensional parameters into the governing equation (Eq.~(\ref{galerkin_eq1})) to facilitate the analysis:
\begin{eqnarray}
& &y=\frac{Y}{H},~ x=\frac{X}{L},~ t=\frac{T}{\alpha},~ \gamma={\frac{EAH^2}{2T_0L^2}}, ~\alpha=\sqrt{\frac{\rho L^2}{T_0}}, \nonumber\\ & & c=\frac{C \alpha}{\rho},~\text{and}~ d(x)=\frac{H-D(x)}{H}.
\label{galerkin_eq3}
\end{eqnarray}
Upon substituting the non-dimensional parameters (Eq.~(\ref{galerkin_eq3})), the governing equation of motion (Eq.~(\ref{galerkin_eq1})) can be written in non-dimensional form as follows:
\begin{equation}
\frac{\partial^2 y}{\partial t^2} - \left(1+\gamma\int_{0}^{L} \left(\frac{\partial y}{\partial x}\right)^2 dx\right) \frac{\partial^2 y}{\partial x^2}+c\frac{\partial y}{\partial t}=0,
\label{galerkin_eq4}
\end{equation}
with the boundary conditions:
\begin{equation}
y(0,t)=y(1,t)=0
\label{galerkin_eq5}
\end{equation}
Now, the solution of Eq.~(\ref{galerkin_eq4}) is assumed as follows:
\begin{equation}
y(x,t)=\sum_{j=1}^{N}\phi_j(x)\eta_j(t) = \boldsymbol{\phi}^\text{T}(x)\boldsymbol{\eta}(t)
\label{galerkin_eq6}
\end{equation}
In Eq.~(\ref{galerkin_eq6}), $\boldsymbol{\phi}(x)=\left[\phi_1(x), \,\ \phi_2(x),...,\phi_N(x)\right]^\text{T}$ and $\boldsymbol{\eta}(t)=\left[\eta_1(t), \,\ \eta_2(t),...,\eta_N(t)\right]^\text{T}$. 
Here, $\phi_j(x)$ are the mass-normalized mode-shapes of the string and $\eta_j(t)$ are the modal coordinates. In this work, $\phi_j(x)$ are chosen to be $\sqrt{2} \sin(j\pi x)$.
On substituting Eq.~(\ref{galerkin_eq6}) in Eq.~(\ref{galerkin_eq4}), pre-multiplying by $\boldsymbol{\phi}(x)$, integrating over the domain $\left[0 ~~1\right]$, and simplifying using orthogonality conditions, the following set of coupled ODEs are obtained:
\begin{equation}
\ddot{\eta_j}(t)+ \left(1 + \gamma \sum_{k=1}^{N} \pi^2 k^2 \eta_k(t)^2\right) \omega_j^2 \eta_j(t)+c\dot{\eta}_j(t)=0
\label{galerkin_eq7}
\end{equation}
Here, $\omega_j=j\pi$ are the natural frequencies of the string. Equation~(\ref{galerkin_eq7}) can now be written in matrix form as
\begin{equation}
\bar{\boldsymbol{M}}\ddot{\boldsymbol{\eta}}(t)+\bar{\boldsymbol{K}}\left(\eta\left(t\right)\right)\boldsymbol{\eta}(t)+\bar{\boldsymbol{C}}\dot{\boldsymbol{\eta}}(t)=0,
\label{galerkin_eq8}
\end{equation}
where $\bar{\boldsymbol{M}}$ is an identity matrix and $\bar{\boldsymbol{K}}(\eta(t))$ is a diagonal matrix with the diagonal elements $\bar{K}_{jj}=\left(1 + \gamma \sum_{k=1}^{N} \pi^2 k^2 \eta_k(t)^2\right) \omega_j^2$ and $\bar{\boldsymbol{C}}=c\bar{\boldsymbol{M}}$. 
Here, Eq.~(\ref{galerkin_eq8}) is in modal coordinates and it is difficult to implement impact conditions in the modal system. Therefore, we discretize the string into a set of physical coordinates (see Fig.~\ref{1d_grid}) and introduce the following notation:
\begin{equation}
y(x_i,t)=p_i(t)
\label{galerkin_eq9}
\end{equation}

Using Eq.~(\ref{galerkin_eq9}), Eq.~(\ref{galerkin_eq6}) can be written as
\begin{equation}
p_i(t)=\boldsymbol{\phi}^\text{T}(x_i)\boldsymbol{\eta}(t)
\label{galerkin_eq10}
\end{equation}
Equation~(\ref{galerkin_eq10}) can be written in matrix form as follows:
\begin{equation}
\boldsymbol{p}(t)=\boldsymbol{\Phi}^\text{T}\boldsymbol{\eta}(t),
\label{galerkin_eq11}
\end{equation}
where $\boldsymbol{p}(t)=\left[p_1(t),p_2(t),...,p_N(t)\right]^\text{T}$ and $\boldsymbol{\Phi}=\left[\boldsymbol{\phi}(x_1),\boldsymbol{\phi}(x_2),...,\boldsymbol{\phi}(x_N)\right]^\text{T}$. Now, Eq.~(\ref{galerkin_eq11}) can be rewritten as
\begin{equation}
\boldsymbol{\eta}(t)=\left(\boldsymbol{\Phi}^\text{T}\right)^{-1}\boldsymbol{p}(t)
\label{galerkin_eq12}
\end{equation}
Upon substituting Eq.~(\ref{galerkin_eq12}) into Eq.~(\ref{galerkin_eq8}) and pre-multiplying by $\boldsymbol{\Phi}^\text{T}$, Eq.~(\ref{galerkin_eq8}) can be written in terms of the physical coordinates as follows:
\begin{equation}
{\boldsymbol{M}}\ddot{\boldsymbol{p}}(t)+{\boldsymbol{K}}\left(\boldsymbol{p}\left(t\right)\right)\boldsymbol{p}(t){+\boldsymbol{C}}\dot{\boldsymbol{p}}(t)=0,
\label{galerkin_eq14}
\end{equation}
where ${\boldsymbol{M}}=\boldsymbol{\Phi}^\text{T}\bar{\boldsymbol{M}}\left(\boldsymbol{\Phi}^\text{T}\right)^{-1}$,  ${\boldsymbol{K}}\left(\boldsymbol{p}\left(t\right)\right)=\boldsymbol{\Phi}^\text{T}\bar{\boldsymbol{K}}\left(\left(\boldsymbol{\Phi}^\text{T}\right)^{-1}\boldsymbol{p}\left(t\right)\right)\left(\boldsymbol{\Phi}^\text{T}\right)^{-1}$ and ${\boldsymbol{C}}=\boldsymbol{\Phi}^\text{T}\bar{\boldsymbol{C}}\left(\boldsymbol{\Phi}^\text{T}\right)^{-1}$.
Now, by applying Ivanov's transformation as discussed in Sec.~\ref{mainder}, to Eq.~(\ref{galerkin_eq14}), the motion of the string impacting a distributed rigid obstacle can be obtained.

\section{Results and discussion}
\label{results}
\begin{figure}[htpb!]
\begin{center}
\includegraphics[width=0.48\textwidth]{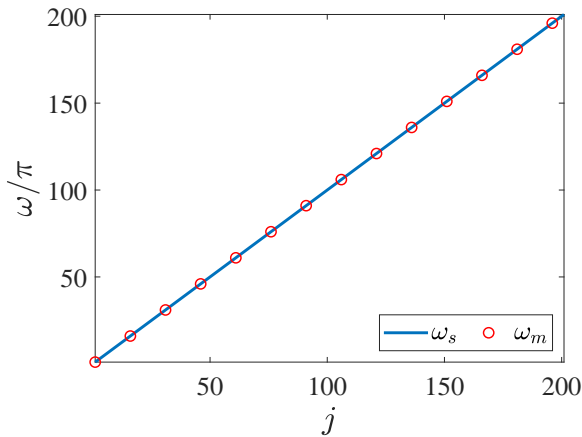}
\caption{Natural frequencies of a string in modal coordinates ${\omega}_m$ (Eq.~(\ref{galerkin_eq8})) and in physical coordinates ${\omega}_s$ (Eq.~(\ref{galerkin_eq14})).}
\label{omegaVsN}
\end{center}
\end{figure}

We now validate the proposed approach by simulating the vibro-impact motion of a string (see Fig.~\ref{Figure7}) against a distributed obstacle using the Galerkin--Ivanov transformation. The analysis is performed by solving Eq.~(\ref{galerkin_eq14}) subjected to the impact constraint using Ivanov's method and penalty method. For the analysis, the non-dimensional parameter $\gamma$ is chosen to be $1$ and $c$ is chosen to be $0$ unless otherwise specified.
The validation is performed by comparing the solutions obtained from Ivanov's method with the results from the penalty approach. The obstacle in the penalty method is modeled as a foundation of unilateral linear springs, each having a spring constant of $k_p=1\times10^{8}$. All numerical simulations have been carried out in MATLAB using 
a fixed-time-step-size, fourth-order Runge--Kutta integrator with $\Delta t = 1 \times 10^{-4}$.
The penalty stiffness for the string and the integration time-step size have been selected to get a good match between the responses from Ivanov's method and penalty method. 
It should be noted that Ivanov's method is a limiting case of the penalty approach (with infinite penalty stiffness) and it is to be expected that the results from these methods will match only for large values of penalty stiffness. 
To demonstrate the reliability of the proposed modal-to-physical transformation, we plotted the first $200$ natural frequencies of the linearized system before and after the transformation in Fig.~\ref{omegaVsN}.
The red circles in Fig.~\ref{omegaVsN} represent the natural frequencies of the linearized system in modal coordinates (${\omega}_m$) (Eq.~(\ref{galerkin_eq8})) and the blue line represents the frequencies in physical coordinates (${\omega}_s$) (Eq.~(\ref{galerkin_eq14})). This clearly demonstrates the accuracy of the transformation in conserving all the natural frequencies of the system exactly.

We now discuss the results for the string impact problem. In the present analysis, we considered $N=201$ modes/grid-points in the Galerkin approximation,for which an accurate match is obtained between Ivanov's method and the penalty approach. As a first example, we simulate the vibro-impacting motion of a string against a flat obstacle.
The considered obstacle is defined as $d(x)=0.025, ~\frac{1}{3}\le x\le \frac{2}{3}$. The initial displacement of the string is considered to be $0.05\times\sin(\pi x)$, and the initial velocity is considered to be zero. On following Ivanov's approach (discussed in Sec.~\ref{mainder}), the solution for the problem is obtained.

\begin{figure}[htpb!]
\begin{center}
\subfigure[]{\includegraphics[width=0.48\textwidth]{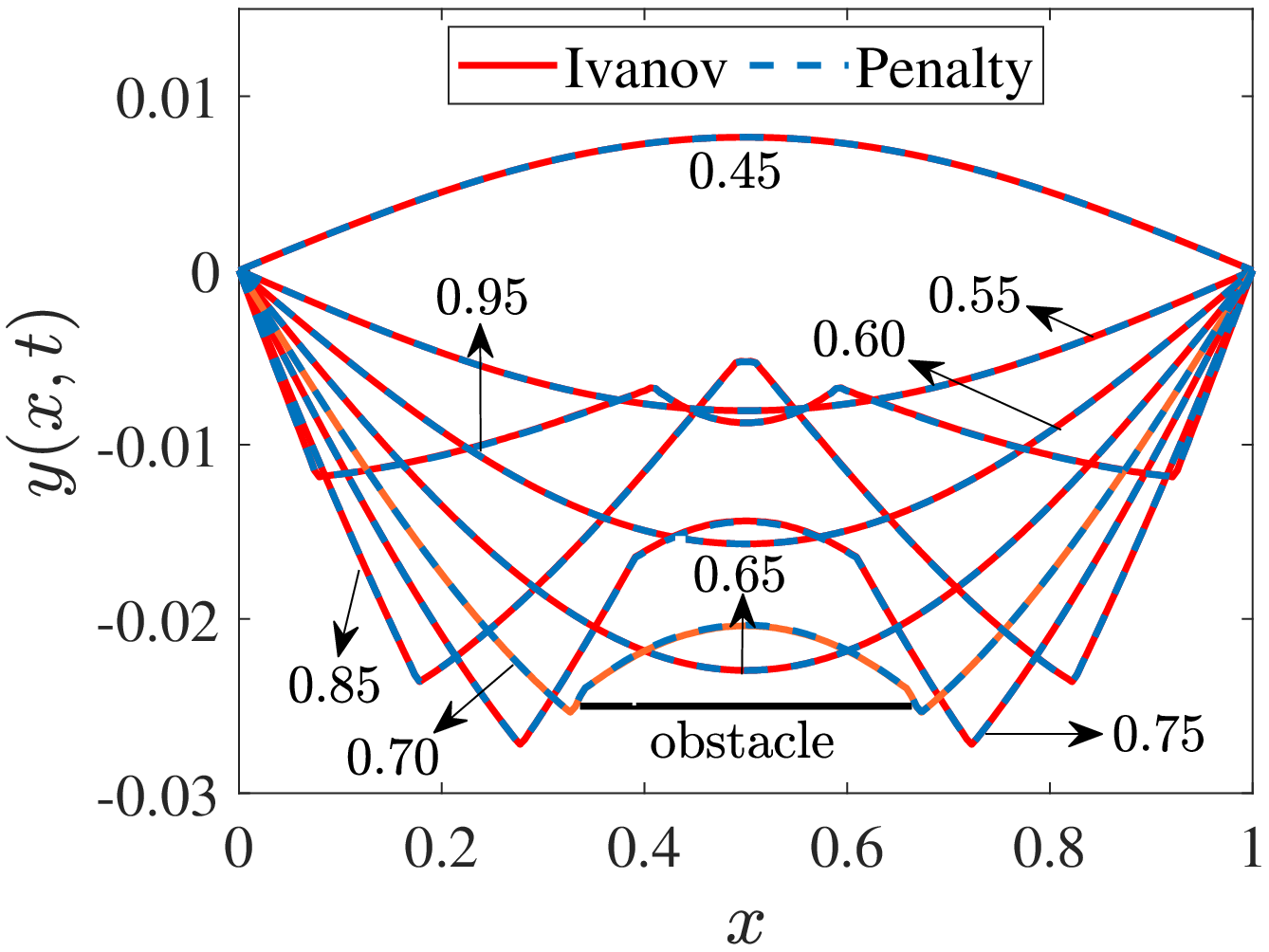}\label{Figure8a}}
\subfigure[]{\includegraphics[width=0.48\textwidth]{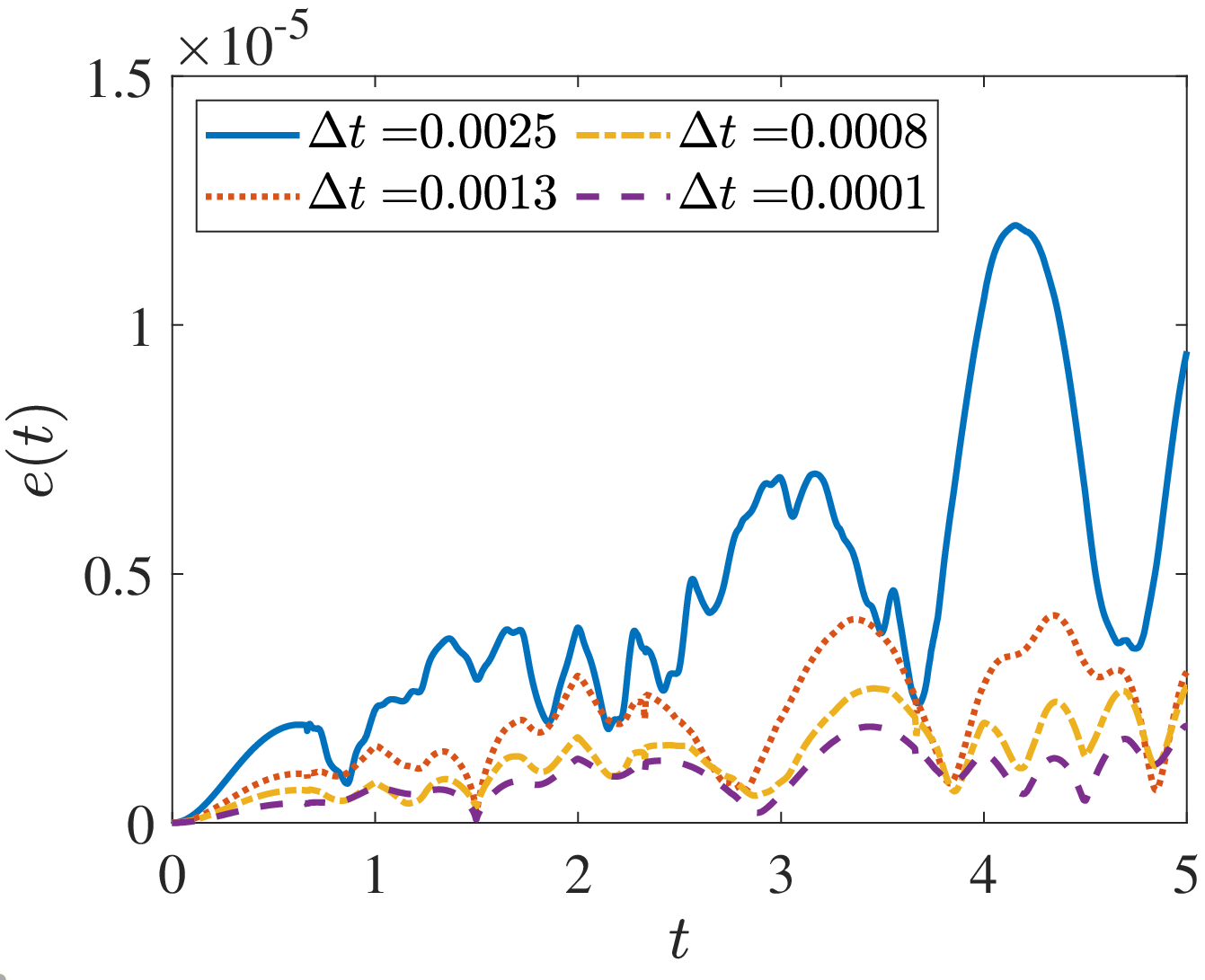}\label{Figure8b}}
\caption{(a) Snapshots of the string motion impacting a flat obstacle at different time instants ($t$) obtained with Ivanov's transformation (solid red line) and the penalty approach (dashed blue line). (b) Mean squared error ($e(t)$) in the string deflection with respect to time ($t$) between the solutions obtained with Ivanov's approach and the penalty method for different integration time-step sizes ($\Delta t$).}
\label{Figure8}
\end{center}
\end{figure}
\begin{figure}[htpb!]
\centering
\subfigure{\includegraphics[width=0.48\textwidth]{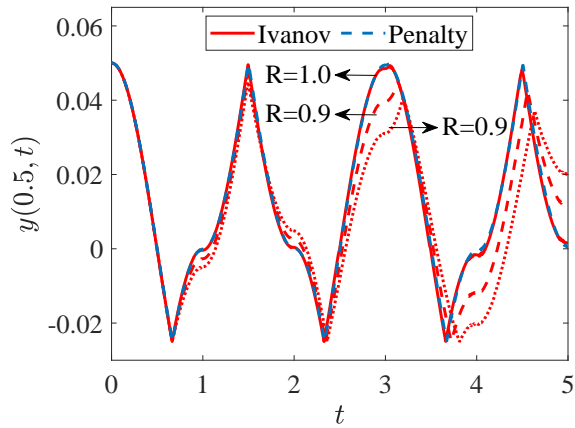}\label{Figure9a}}
\caption{Midpoint displacement of the string impacting a flat obstacle obtained with Ivanov's transformation: $R=1$ (solid red line), $R=0.9$ (dashed red line), and $R=0.8$ (dotted red line). The dashed blue lines show the results obtained using the penalty approach.}
\label{Figure9}
\end{figure}

In Fig.~\ref{Figure8a}, we have shown the shape of the string obtained from numerical simulations at different time instants for $R=1$ (CoR) obtained using both the Ivanov and penalty approaches. A clear match between the results obtained from both the methods has been observed in Fig.~\ref{Figure8a}. Furthermore, we have also implemented the Ivanov approach using an explicit numerical integrator for different time-step sizes ($\Delta t=~0.0025,~ 0.0013,~ 0.0008$, and $0.0001$). The MSE ($e(t)$) between the solutions obtained from Ivanov and penalty approaches with respect to time has been evaluated and shown in Fig.~\ref{Figure8b}. The solution from the penalty approach has been determined for the penalty spring stiffness of $k_p=10^8$ and is used to compare the Ivanov solutions. From Fig.~\ref{Figure8b}, it is clear that the error in the Ivanov solution reduces with decrease in the integration time-step size and converges to the penalty solution. It is also observed that the initial MSE between the solutions is of an order less than $10^{-5}$ and increases gradually with respect to time. This happens because of the finite penalty stiffness, due to which the duration of contact is also finite. When the simulation is run for a long time, these accumulated contact times will introduce a phase shift in the solution when compared to Ivanov's solution. Moreover, for simulating a rigid collision in penalty approach, one has to use a large value for the penalty stiffness, which puts a hard restriction on the numerical integrator time-step size. However, in Ivanov's approach, the nonsmooth transformation simulates the case of infinite penalty stiffness, but without making the equations stiff. Therefore, the comparison between the Ivanov and penalty methods will be valid only for small simulation times (few impacts). In addition, the unilateral constraints are exactly satisfied in Ivanov's approach, unlike the penalty method. We can only compare results obtained from Ivanov's method for the case of $R=1$ with those obtained from the penalty approach. This is because there are no models that relate $R$ with contact dissipation in the penalty approach for discretized continuous systems.

\begin{figure}[htpb!]
\centering
\subfigure[]{\includegraphics[width=0.48\textwidth]{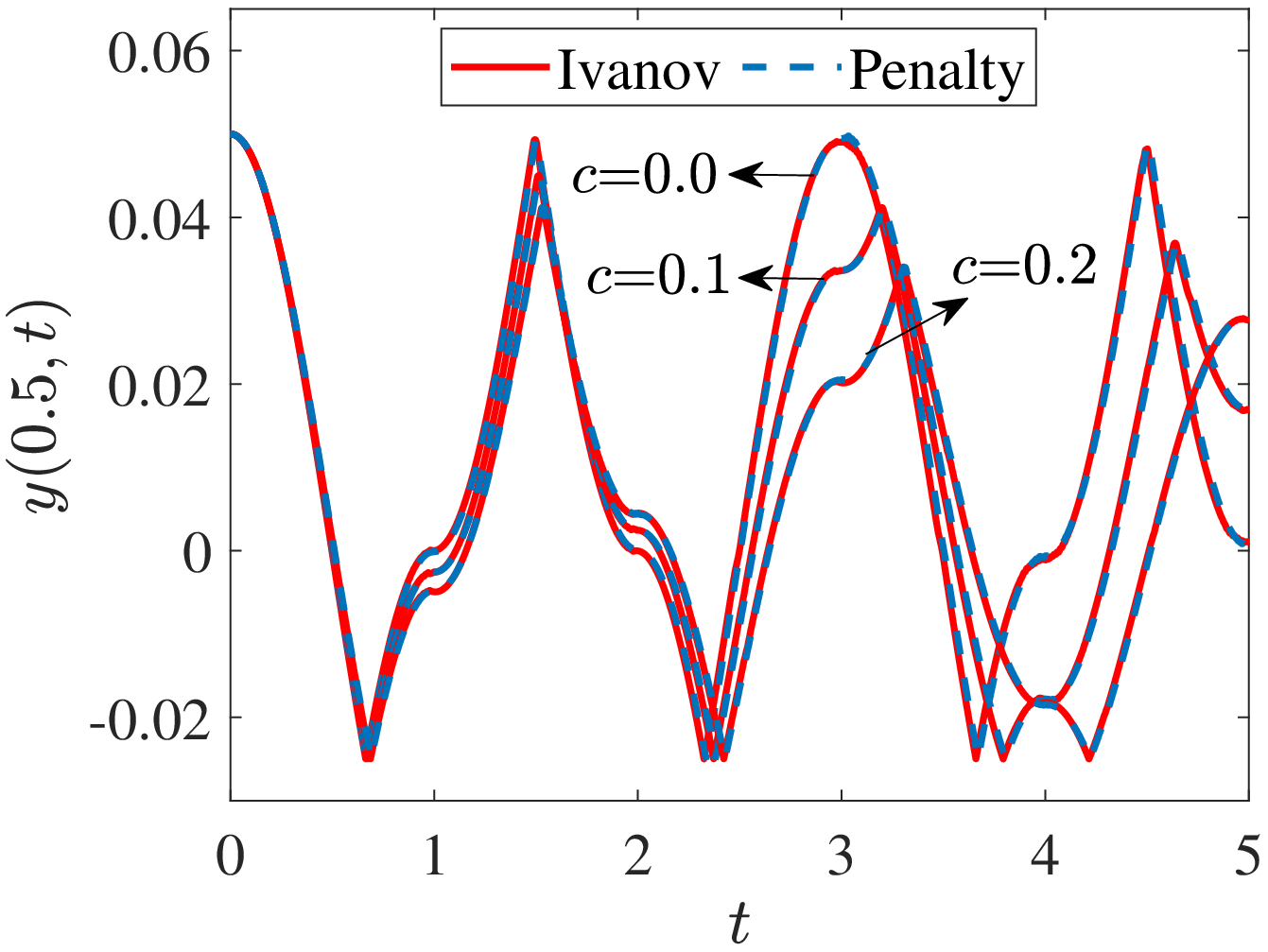}\label{Figure12a}}
\subfigure[]{\includegraphics[width=0.48\textwidth]{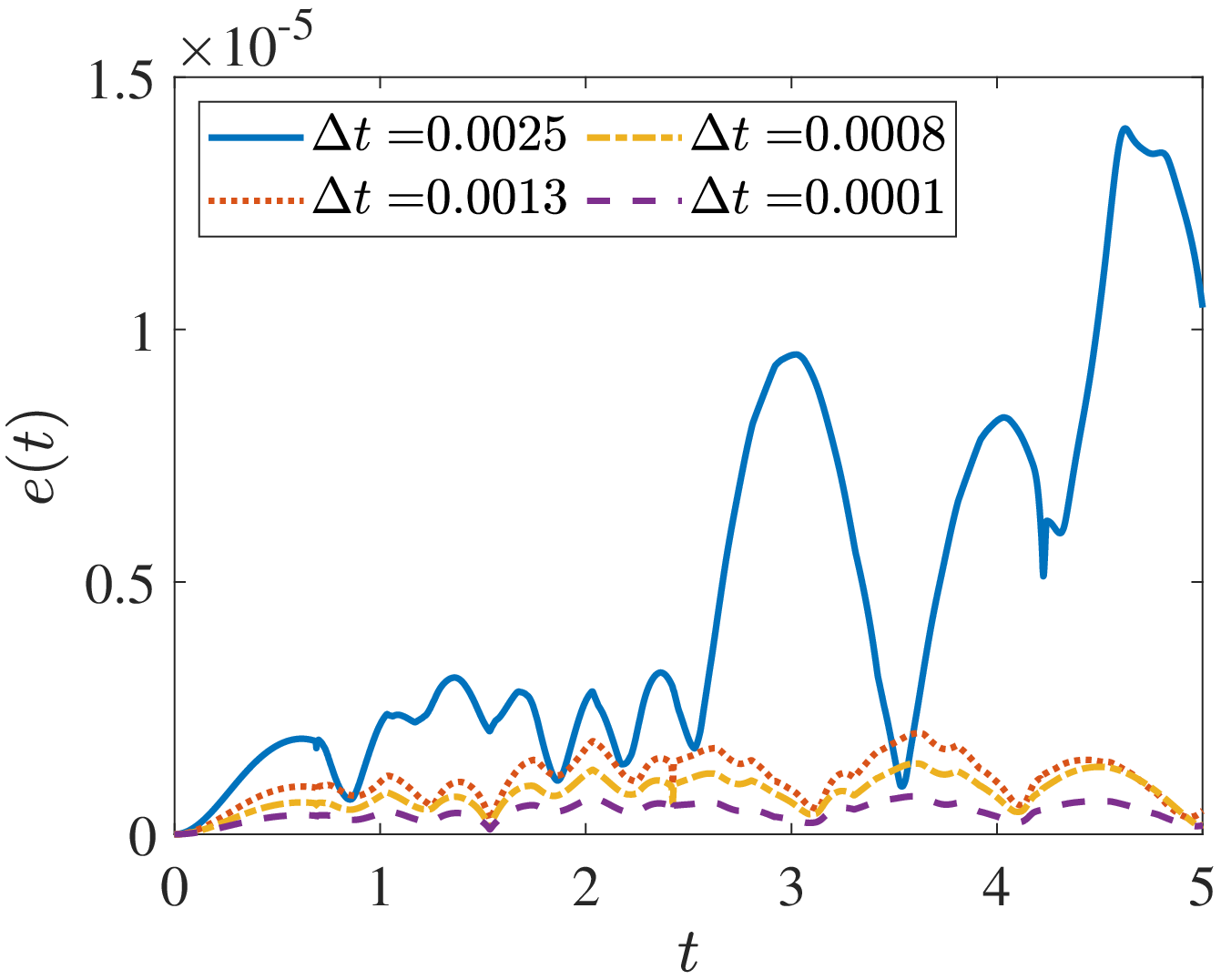}\label{Figure12b}}
\caption{(a) Midpoint displacement of the string impacting a flat obstacle obtained with Ivanov's transformation for $c=0.0,~0.1,$ and $0.2$. The dashed blue lines show the results obtained using the penalty approach. (b) Error between Ivanov and penalty solutions for different integration time-step sizes ($\Delta t$) for $c=0.2$ and $R=1$.}
\label{Figure12}
\end{figure}

In Fig.~\ref{Figure9}, we show the midpoint deflection of the string for three values of CoR ($R= 1.0$, $0.9$, and $0.8$) for $c=0$. Again, we can see a good agreement between the solutions for the penalty method (dashed blue line) and Ivanov's method for $R=1$ (solid red line). It is also observed from Fig.~\ref{Figure9} that the rebounding displacement decreases for a smaller $R$ because of the loss of energy with every impact for $R<1$. We also present the midpoint displacement of the string impacting a flat obstacle for different damping values ($c=0,~0.1$, and $0.2$) in Fig.~\ref{Figure12a}. The corresponding penalty solutions have also been shown in Fig.~\ref{Figure12b} with a dashed blue line and an accurate match between both the solutions has been observed. In addition, the MSE between the Ivanov and penalty solutions for different integration time-step sizes ($\Delta t=0.0025,~0.0013,~0.0008$ and $0.0001$) for $c=0.2$ and $R=1$ has been presented in Fig.~\ref{Figure12b}. It is clearly observed from Fig.~\ref{Figure12b} that the MSE between both solutions is of the order less than $1.5\times 10^{-5}$ for all $\Delta t$. Also the Ivanov solution seems to converge to the penalty solution for decreasing $\Delta t$.

\begin{figure}[htpb!]
\begin{center}
\subfigure[]{\includegraphics[width=0.48\textwidth]{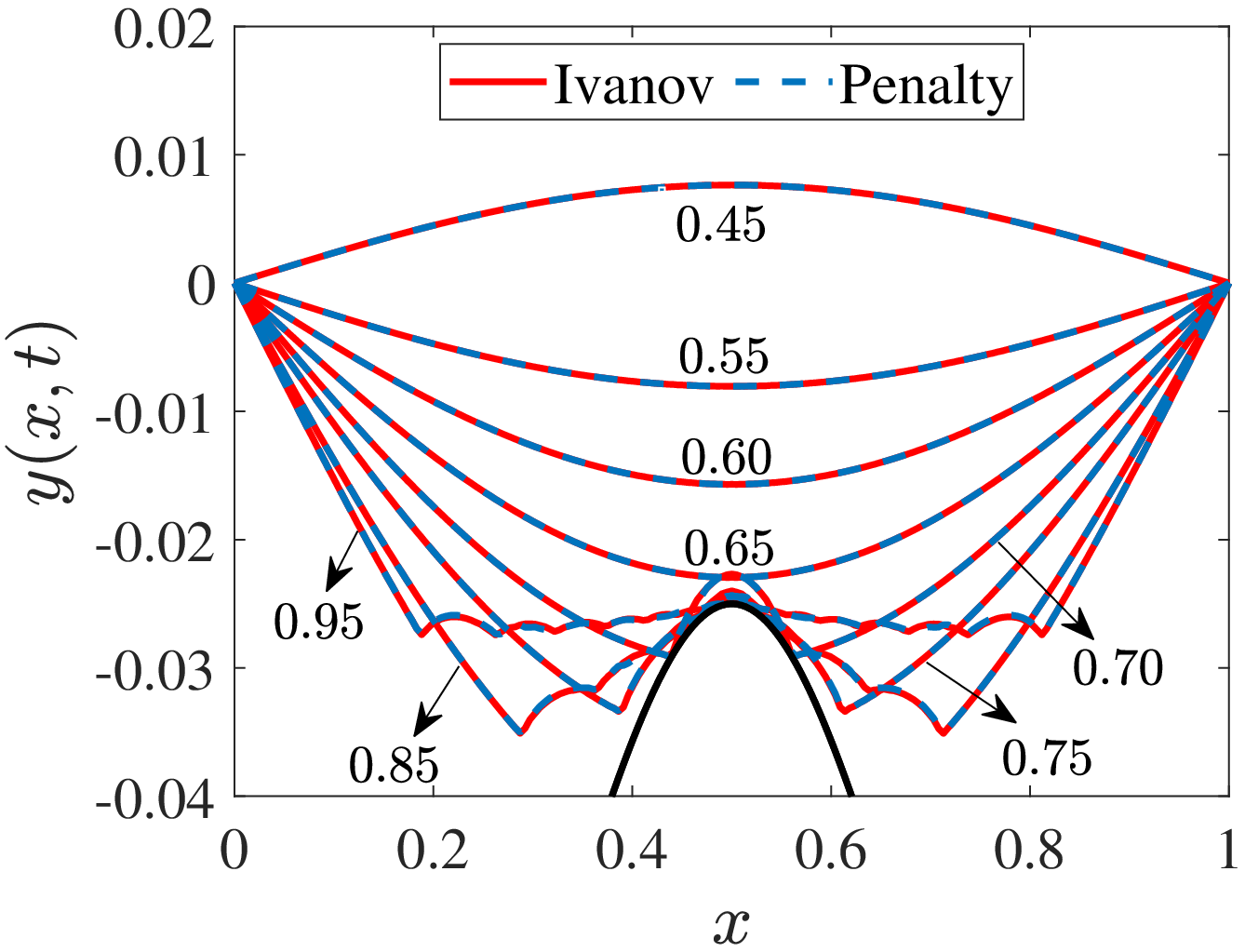}\label{Figure10a}}
\subfigure[]{\includegraphics[width=0.48\textwidth]{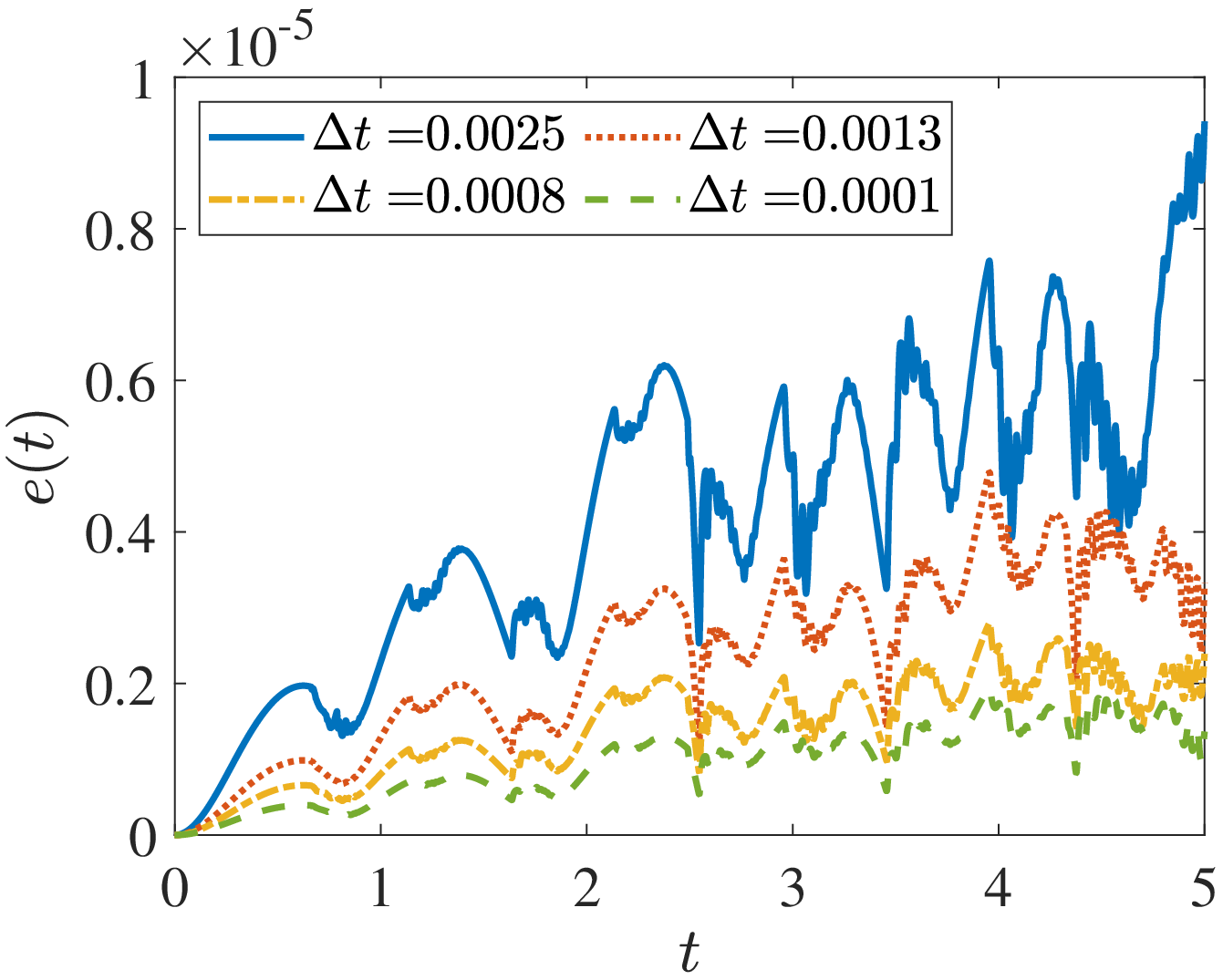}\label{Figure10b}}
\caption{(a) Snapshots of the string motion impacting a sinusoidal obstacle obtained with Ivanov's transformation (solid red line) and the penalty approach (dashed blue line). (b) Error in the string deflection between the solutions obtained from Ivanov's approach and the penalty method.}
\label{Figure10}
\end{center}
 \end{figure}
\begin{figure}[htpb!]
\centering
\subfigure{\includegraphics[width=0.48\textwidth]{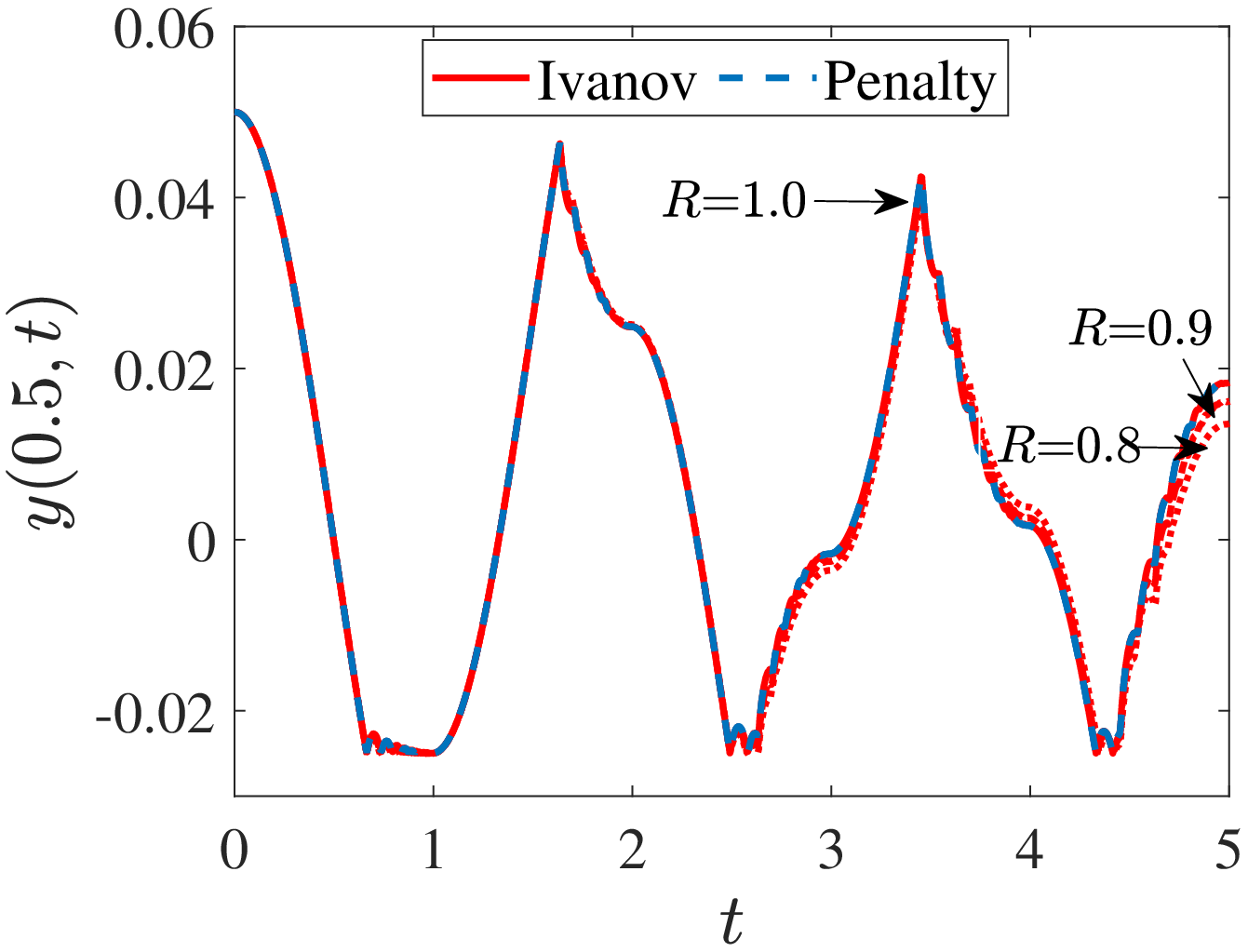}\label{Figure11a}}
\caption{Midpoint displacement of the string impacting a sinusoidal obstacle obtained with Ivanov's transformation: $R=1$ (solid red line), $R=0.9$ (dashed red line), and $R=0.8$ (dotted red line). The dashed blue lines show the results obtained using the penalty approach.}
\label{Figure11}
\end{figure}

\begin{figure}[htpb!]
\centering
\subfigure[]{\includegraphics[width=0.48\textwidth]{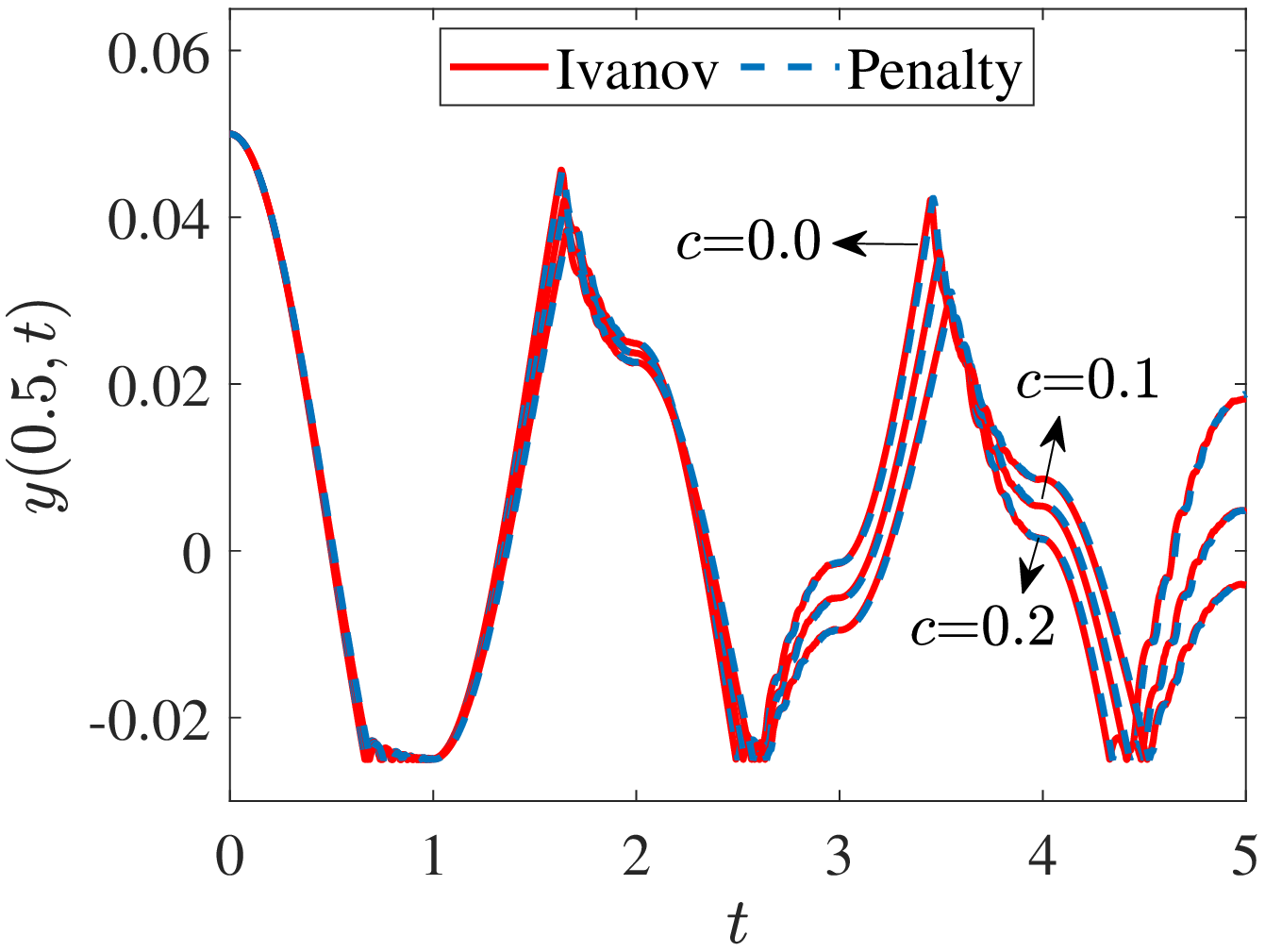}\label{Figure13a}}
\subfigure[]{\includegraphics[width=0.48\textwidth]{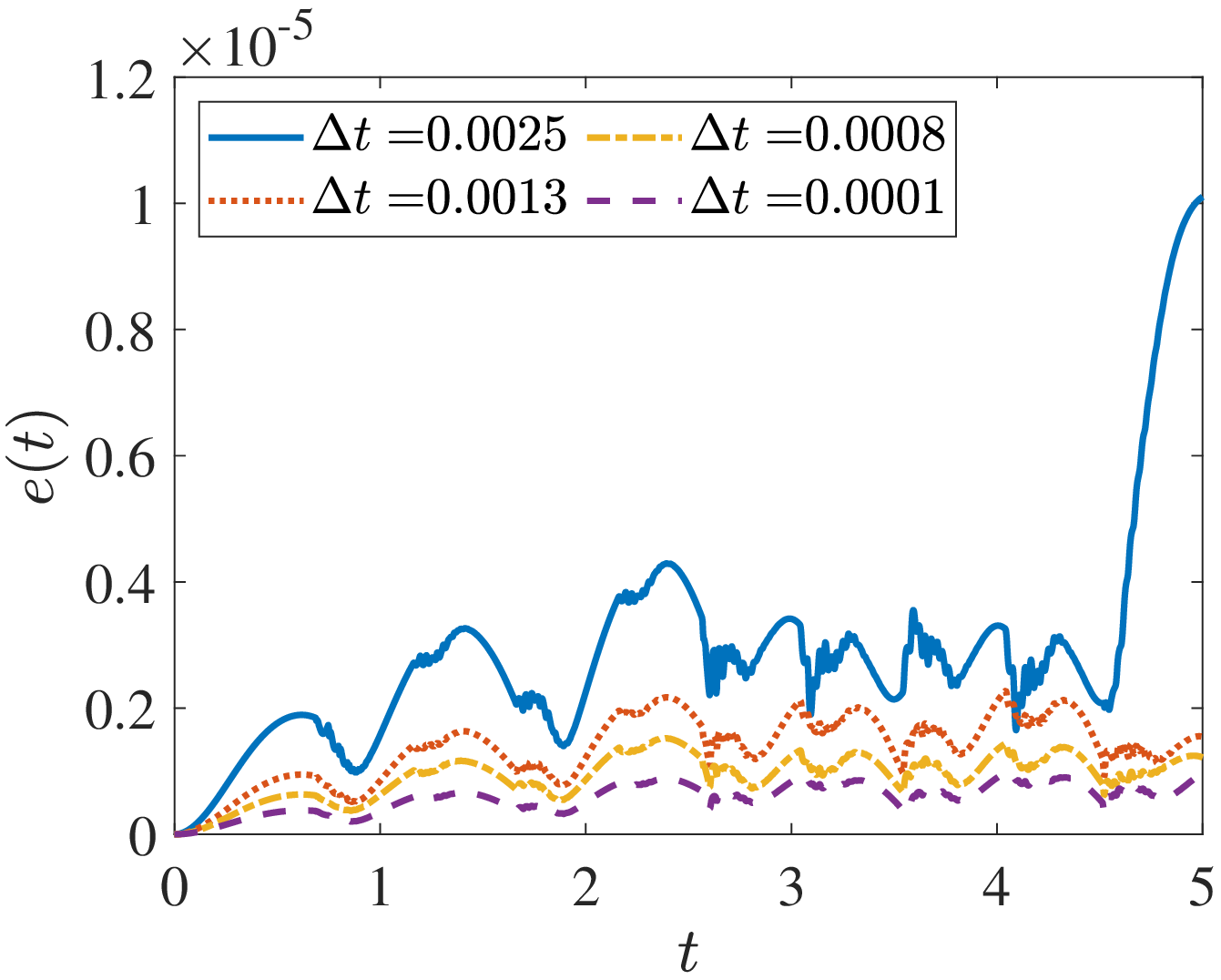}\label{Figure13b}}
\caption{(a) Midpoint displacement of the string impacting a sinusoidal obstacle obtained with Ivanov's transformation for $c=0.0,~0.1,$ and $0.2$. The dashed blue lines show the results obtained using the penalty approach. (b) Error between Ivanov and penalty solutions for different integration time-step sizes ($\Delta t$) for $c=0.2$ and $R=1$.}
\label{Figure13}
\end{figure}

For the second example, we consider the motion of a string impacting a sinusoidal obstacle of the form $d(x)=0.05-0.025\times\sin(\pi (x-\frac{1}{3})), ~\frac{1}{3}\le x\le \frac{4}{3}$.
The initial displacement of the string is considered to be $0.05 \times \sin(\pi x)$, and the initial velocity is taken to be zero. Figures~\ref{Figure10}, \ref{Figure11}, and \ref{Figure13} are similar to Figs.~\ref{Figure8}, \ref{Figure9}, and \ref{Figure12}, respectively, except that the results are presented for the case of a sinusoidal obstacle. It is clearly observed from Figs.~\ref{Figure10}, \ref{Figure11} and \ref{Figure13} that the solutions from Ivanov's method and the penalty approach match closely with an MSE of order less than $10^{-5}$. 

From the results shown in Figs.~\ref{Figure8}--\ref{Figure11}, it is clearly demonstrated that Ivanov's transformation can be successfully used to simulate the vibro-impacting motions of a continuous structure with distributed contacts using a Galerkin approximation. Even though we have reported results for a string, the approach can also be used to simulate impacts in membranes, beams, and plate structures.

\section{Conclusions}
\label{concl}
In this paper, the application of Ivanov's method to simulate the nonsmooth motion of a continuous structure against a rigid distributed obstacle using Galerkin approximation is demonstrated. A closed-form expression for the nonlinear MDOF system in Ivanov's transformed coordinates is reported. The transformed differential equations are without any constraints and automatically satisfy the unilateral and velocity jump conditions at impact. The technique is validated by solving three problems: (i) a spring-mass-damper system, (ii) the motion of a string against a flat obstacle, and (iii) the motion of a string against a sinusoidal obstacle. The governing partial differential equations of motion of the string are approximated using the Galerkin approach. We first convert the ODEs obtained from the Galerkin approximation to the physical coordinates using a modal-physical transformation and then apply Ivanov's approach. 
The solutions that are obtained from Ivanov's method are then compared with the results from the penalty method.
The equations of motion obtained from Ivanov's method are not stiff, in contrast to those obtained from the penalty method. Therefore, a larger time-step size can be used for the numerical integration of the equations of motion obtained from Ivanov's method. The event detection problem for large MDOF systems in the presence of multiple impact constraints is a difficult and challenging task and is completely eliminated by the application of Ivanov's method.

\section{Acknowledgments} \label{acknowledgments}
CPV gratefully acknowledges the Department of Science and Technology for funding this research through Inspire fellowship (DST/INSPIRE/04/2014/000972).

\bibliographystyle{ieeetr}

\end{document}